\documentclass[10pt,draft]{amsart}
\usepackage{amsmath,amssymb,amsthm}
\begin{document}
\newtheorem{lem}{Lemma}[section]
\newtheorem{prop}{Proposition}[section]
\newtheorem{cor}{Corollary}[section]
\numberwithin{equation}{section}
\newtheorem{thm}{Theorem}[section]
\theoremstyle{remark}
\newtheorem{example}{Example}[section]
\newtheorem*{ack}{Acknowledgment}
\theoremstyle{definition}
\newtheorem{definition}{Definition}[section]
\theoremstyle{remark}
\newtheorem*{notation}{Notation}
\theoremstyle{remark}
\newtheorem{remark}{Remark}[section]
\newenvironment{Abstract}
{\begin{center}\textbf{\footnotesize{Abstract}}%
\end{center} \begin{quote}\begin{footnotesize}}
{\end{footnotesize}\end{quote}\bigskip}
\newenvironment{nome}
{\begin{center}\textbf{{}}%
\end{center} \begin{quote}\end{quote}\bigskip}

\newcommand{\triple}[1]{{|\!|\!|#1|\!|\!|}}
\newcommand{\xx}{\langle x\rangle}
\newcommand{\ep}{\varepsilon}
\newcommand{\al}{\alpha}
\newcommand{\be}{\beta}
\newcommand{\de}{\partial}
\newcommand{\la}{\lambda}
\newcommand{\La}{\Lambda}
\newcommand{\ga}{\gamma}
\newcommand{\del}{\delta}
\newcommand{\Del}{\Delta}
\newcommand{\sig}{\sigma}
\newcommand{\ome}{\Omega^n}
\newcommand{\Ome}{\Omega^n}
\newcommand{\C}{{\mathbb C}}
\newcommand{\N}{{\mathbb N}}
\newcommand{\Z}{{\mathbb Z}}
\newcommand{\R}{{\mathbb R}}
\newcommand{\T}{{\mathbb T}}
\newcommand{\Rn}{{\mathbb R}^{n}}
\newcommand{\Rnu}{{\mathbb R}^{n+1}_{+}}
\newcommand{\Cn}{{\mathbb C}^{n}}
\newcommand{\spt}{\,\mathrm{supp}\,}
\newcommand{\Lin}{\mathcal{L}}
\newcommand{\SSS}{\mathcal{S}}
\newcommand{\F}{\mathcal{F}}
\newcommand{\xxi}{\langle\xi\rangle}
\newcommand{\eei}{\langle\eta\rangle}
\newcommand{\xei}{\langle\xi-\eta\rangle}
\newcommand{\yy}{\langle y\rangle}
\newcommand{\dint}{\int\!\!\int}
\newcommand{\hatp}{\widehat\psi}
\renewcommand{\Re}{\;\mathrm{Re}\;}
\renewcommand{\Im}{\;\mathrm{Im}\;}

\title[The NLS ground states on product spaces]%
{The Nonlinear Schr\"odinger  
equation ground states on product spaces}

\author[Susanna Terracini]{Susanna Terracini}
\author[Nikolay Tzvetkov]{Nikolay~Tzvetkov}
\author[Nicola Visciglia]{Nicola Visciglia}
\address{Dipartimento di Matematica e Applicazioni, Universit\`a  di Milano Bicocca, Via Cozzi 53,
20125 Milano, Italy}
\email{susanna.terracini@unimib.it }
\address{D\'epartement de Math\'ematiques, Universit\'e de Cergy-Pontoise, 2, avenue Adolphe Chauvin, 95302 Cergy-Pontoise  
Cedex, France and Institut Universitaire de France}
\email{nikolay.tzvetkov@u-cergy.fr}
\address{Universit\`a Degli Studi di Pisa Dipartimento di Matematica "L. Tonelli"
Largo Bruno Pontecorvo 5 I - 56127 Pisa. Italy}
\email{ viscigli@dm.unipi.it}
\maketitle
\date{}
\begin{abstract}
We study the nature of the Nonlinear Schr\"odinger equation ground states on the product spaces $\R^n\times M^k$, where $M^k$ is a compact Riemannian manifold. 
We prove that for small $L^2$ masses the ground states coincide with the corresponding $\R^n$ ground states. We also prove that above a critical mass the ground states have nontrivial $M^k$ dependence. Finally, we address the Cauchy problem issue which transform the variational analysis to dynamical stability results.
\end{abstract}
\section{Introduction}
Our goal here is to study the nature of the 
Nonlinear Schr\"odinger equation ground states when 
the problem is posed on the product spaces $\R^n\times M^k$, 
where $M^k$ is a compact Riemannian manifold. We thus consider the 
following Cauchy problems
\begin{equation} \label{NLS}
\left\{ \begin{array}{c}
i\partial_t u - \Delta_{x,y} u - u|u|^{\alpha}=0, \hbox{ } 
(t, x,y)\in \R \times \R^{n}_x\times M^k_y
\\
u(0,x,y)= \varphi(x, y)
\end{array} \right.
\end{equation}
where 
$$
\Delta_{x,y}= \sum_{j=1}^{n}\partial_{x_j}^2+ \Delta_y 
$$ 
and $\Delta_y$ is the Laplace-Beltrami operator on $M^k_y$.
Recall that the Laplace-Beltrami operator is defined
in local coordinates as follows:
$$\frac 1{\sqrt {det (g_{i,j}(y))}} \partial_{y_i} 
\sqrt {det (g_{i,j}(y))} g^{i,j}(y) \partial_{y_j}$$
where
$g^{i,j}(y)= (g_{i,j}(y))^{-1}$ and $g_{i,j}(y)$ is the metric tensor.\\
We assume that $0<\alpha<4/(n+k)$ which corresponds to $L^2$ subcritical nonlinearity.
In this paper, we shall study the following two questions: 
\begin{itemize}\item the existence and 
stability of solitary waves for \eqref{NLS};
\item the global well posedness of the Cauchy problem associated to \eqref{NLS}.
\end{itemize}
The equation \eqref{NLS} has two (at least formal) conservation laws, the energy
\begin{equation}\label{L2peri}\mathcal 
E_{n, M^k, \alpha}(u)= \int_{M^k_y} \int_{\R_x^{d}} 
\Big (\frac 12 |\nabla_{x, y} u|^2 
-\frac{1}{2+\alpha} |u|^{2+\alpha} \Big )dxdvol_{M^k_y}
\end{equation}
and the $L^2$ mass,
\begin{equation}\label{constrain}\|u\|_{L^2(\R^n\times M^k)}^2=
\int_{M^k_y} \int_{\R^{n}_x} |u|^2 dxdvol_{M_y^k}
\end{equation}
Here we denote by $dvol_{M_y^k}$ the volume form on $M^k$. Recall that 
in local coordinates it can be written
as 
$\sqrt {det (g_{i,j}(y))} dy$. Moreover 
the $i$-th component (in local coordinates) of the gradient $(\nabla_y u(y))$ 
is $$g^{i,j}(y)\partial_{y_j}u$$
One has the classical Gagliardo-Nirenberg inequality
\begin{equation}\label{GagN}
\|u\|_{L^{2+\alpha}( \R^n\times M^k)}^{2+\alpha}\leq C
\|u\|_{H^{1}( \R^n\times M^k)}^{\theta(\alpha)}
\|u\|_{L^{2}( \R^n\times M^k)}^{2+\alpha-\theta(\alpha)}
\end{equation}
where
$
\theta(\alpha)=(n+k)\alpha/2$. Thus $\theta(\alpha)<2$ under 
our assumption $0<\alpha<4/(n+k)$.
This implies that the conservation laws \eqref{L2peri} and \eqref{constrain} 
imply a control on the $H^1$ norm which excludes a $L^2$ 
self-focusing blow-up and thus
one expects that \eqref{NLS} has a well-defined global dynamics. 
This problem seems quite delicate for a general $M^k$. However 
if we replace $M^k$ with $\R^k$ it is well-known (see \cite{Ts}, \cite{C} 
and the references therein) that \eqref{NLS} 
has a global strong solution for every $L^2(\R^{n+k})$ initial data.   

Our argument to construct stable solutions to \eqref{NLS}
follows the one proposed in \cite{CL}. Hence we shall look at the following 
minimization problems: 
\begin{equation}\label{minta}K_{n, M^k, \alpha}^\rho=
\inf_{\substack{u\in H^1(\R^n\times M^k)\\
\|u\|_{L^2(\R^n \times M^k)}=\rho}}{\mathcal E}_{n, M^k, \alpha}(u)
\end{equation}
and ${\mathcal E}_{n, M^k, \alpha}(u)$ is defined in \eqref{L2peri}.
\noindent In the sequel we shall use the following notation: 
\begin{multline}\label{minim56}
{\mathcal M}_{n, M^k, \alpha}^\rho
=\{v\in H^1(\R^n\times M^k)|
\\
 \|v\|_{L^2(\R^n \times M^k)}=\rho
\hbox{ and } {\mathcal E}_{n,M^k,\alpha}(v)=K^\rho_{n, M^k,\alpha}\}
\end{multline}
The first result we state concerns the compactness of 
minimizing sequences to \eqref{minta}.
\begin{thm}\label{variationalappendix}
Let $M^k$ be a compact manifold and $0<\alpha<4/(n+k)$. 
Then we have the following:
\begin{equation}\label{1}
K^{\rho}_{n, M^k, \alpha}>-\infty \hbox{ and }
{\mathcal M}_{n, M^k, \alpha}^{\rho}\neq \emptyset, \hbox{ } \forall \rho>0;
\end{equation}
\begin{equation}\label{2}
\forall u_j \in H^1(\R^n\times M^k) \hbox{ s.t. } 
\|u_j\|_{L^2(\R^n\times M^k)}=\rho, \lim_{j\rightarrow \infty}
{\mathcal E}_{n, M^k,\alpha}(u_j)= K^{\rho}_{n, M^k, \alpha}
\end{equation}
$$\hbox{ } \exists \hbox{ a subsequence } u_{j_l}
\hbox{ and }\tau_l\in \R_x^{n} \hbox{ s.t. }
u_{j_l} (x+ \tau_l, y) \hbox{ converges in } H^1(\R^{n}\times M^k).$$
\end{thm}
The proof of Theorem \ref{variationalappendix} is based on the concentration compactness principle and it 
will be given in the appendix. Also the following stability theorem
follows from a standard argument, hence its classical proof will be recalled in the appendix.
\begin{thm}\label{stabilitysubordinate}
Let $\rho>0$ be fixed and $n, M^k, \alpha$ as in Theorem \ref{variationalappendix}.
Assume moreover that 
\begin{equation}\label{gwp}\hbox{ the Cauchy problem \eqref{NLS} is globally well posed
for any data $\varphi\in \mathcal U$}\end{equation}$$\hbox{ where $\mathcal U$ is a $H^1(\R^n\times M^k)$-neighborhood of
${\mathcal M}_{n, M^k, \alpha}^\rho$.}$$ Then the set
${\mathcal M}_{n, M^k, \alpha}^\rho$ is orbitally stable, i.e.:
$$\forall \epsilon>0 \hbox{ } \exists \delta=\delta(\epsilon)>0 \hbox{ s.t. }
$$ 
$$\varphi\in \mathcal U, \inf_{v \in
{\mathcal M}_{n, M^k, \alpha}^\rho} \|\varphi - v\|_{H^1(\R^n\times M^k)}<\delta(\epsilon)$$ 
$$ \hbox{ implies } \sup_{t\in \R} \Big (\inf_{v \in
{\mathcal M}_{n, M^k, \alpha}^\rho} 
\|u_\varphi(t) - v\|_{H^1(\R^n\times M^k)}\Big)<\epsilon$$
where $u_\varphi(t,x,y)$ is the unique global solution to \eqref{NLS}.
\end{thm}
Let us emphasize that the stability result stated 
in Theorem \ref{stabilitysubordinate}
has two major defaults: the first one is that we don't have an explicit description of the minimizers
${\mathcal M}_{n, M^k, \alpha}^\rho$;
the second one is that it is subordinated to \eqref{gwp}, i.e. the 
global well posedness of the Cauchy problem \eqref{NLS}.
The main contributions of this paper concern a partial understanding of the aforementioned questions.
\\
Notice that (see \cite{C}) a special family of solutions to \eqref{NLS} is given by
$$u(t, x, y)=e^{-i\omega t} u_{n,\omega,\alpha}(x)$$ 
where $\omega>0$ and $u_{n,\omega, \alpha}(x)$ is defined as the unique radial solution to: 
\begin{equation}\label{uomegaintro}
-\Delta_x u_{n,\omega,\alpha} + \omega u_{n,\omega,\alpha}= u_{n,\omega,\alpha}
|u_{n,\omega,\alpha}|^{\alpha}
\end{equation}$$
 \hbox{ }
 u_{n,\omega, \alpha}\in H^1(\R^{n}_x), \hbox{ }
u_{n,\omega, \alpha}(x)>0, \hbox{  }  x\in \R_x^{n}
$$
Next, we set
\begin{equation}\label{Nomega}
{\mathcal N}_{n,\omega,\alpha}=\{e^{i\theta} u_{n,\omega,\alpha}(x+\tau)|
\tau\in \R^{n}, \theta \in \R\}
\end{equation}
Notice that there is a natural embedding  $H^1(\R^{n}_x)\subset
H^1(\R^n_x\times M^k_y)$. In fact every function in $H^1(\R^{n}_x)$
can be extended in a trivial way w.r.t. the $y$ variable on $\R^n_x\times M^k_y$ and 
this extension will belong to
$H^1(\R^n\times M^k)$.
In particular since now on the set ${\mathcal N}_{n,\omega, \alpha}$ defined in \eqref{Nomega},
will be considered without any further comment in a twofold way: 
as a subset of $H^1(\R^{n}_x)$ and $H^1(\R^{n}_x\times M^k_y)$.
By a rescaling argument
one can prove that the function
$$(0, \infty) \ni \omega \rightarrow 
\|u_{n,\omega,\alpha}\|_{L^2(\R^{n}_x)}^2
\in (0,\infty)$$
is strictly increasing for any $0<\alpha<\frac 4n$ and
$$\lim_{\omega \rightarrow \infty }
\|u_{n,\omega, \alpha}\|_{L^2(\R^{n}_x)}= \infty \hbox{ and }
\lim_{\omega \rightarrow 0} \|u_{n,\omega, \alpha}\|_{L^2(\R^{n}_x)}=0$$
As a consequence for any fixed $0<\alpha<\frac 4n$ we have:
\begin{equation}\label{remtriv}\forall \rho>0 \hbox{ } 
\exists!  \hbox{ } 
\omega(\rho)>0 \hbox{ s.t. } 
\|u_{n,\omega(\rho), \alpha}\|_{L^2(\R^{n}_x)}=\rho
\end{equation}
In next theorem the set 
${\mathcal N}_{n,\omega,\alpha}$ is the one defined in \eqref{Nomega}
and ${\mathcal M}_{n, M^k, \alpha}^\rho$ is defined in \eqref{minim56}.
\begin{thm}\label{variational}
Let $n, M^k, \alpha$ as in Theorem \ref{stabilitysubordinate}.
There exists $\rho^*\in (0, \infty)$ such that:
\begin{equation}\label{rigi}
{\mathcal M}_{n, M^k, \alpha}^{\rho}= 
{\mathcal N}_{n,\omega(\rho/\sqrt {vol(M^k)}),\alpha}, 
\hbox{ } \forall \rho<\rho^*\end{equation}
and
\begin{equation}\label{norigi}
{\mathcal M}_{n, M^k, \alpha}^{\rho}\cap 
{\mathcal N}_{n,\omega(\rho/\sqrt {vol(M^k)}),\alpha}=\emptyset, 
\hbox{ } \forall \rho>\rho^*\end{equation}
where $\omega(\rho /\sqrt {vol(M^k)})$ is 
uniquely defined in \eqref{remtriv}.
In particular for $\rho>\rho^*$ the elements of ${\mathcal M}_{n, M^k, \alpha}^{\rho}$ depend in a nontrivial way on the $M^k$ variable.
\end{thm}
By the approach of Weinstein \cite{W} one may expect that ${\mathcal N}_{n,\omega,\alpha}$ 
is stable under \eqref{NLS} for $\alpha<4/n$ and $\omega$ small enough, 
see \cite{RT} for a recent related work. 
It should however be pointed out that in such a stability result one would not get the variational description of 
${\mathcal N}_{n,\omega,\alpha}$  
as is the case in Theorem~\ref{variational} ($\alpha<4/(n+k)$).
We underline that by combining Theorem \ref{stabilitysubordinate}
and Theorem \ref{variational} we get a stable set for large values of the mass $\rho$,
and in general it is independent of the solitary solitary waves 
associated to NLS in $\R^n$.
\\
Next we shall focus on the question of the global well-posedness of the Cauchy problem
associated to \eqref{NLS} in the particular case $n\geq 1$, $k=1$.
For every $n>1$ we fix the numbers
$$p:=p(n, \alpha)=\frac{4(2+\alpha)}{n\alpha} \hbox{ and }
q:=q(n, \alpha)=2+\alpha$$ and for every $T>0$ we define
the localized norms:
\begin{equation}\label{XT}\|u(t,x,y)\|_{X_{T}}\equiv\|u(t,x,y)
\|_{L^p((-T, T); L^q(\R^{n}_x;H^1(M^1_y))}\end{equation}
and 
\begin{equation}\label{YT}\|u(t,x,y)\|_{Y_{T}}\equiv\|\nabla_x 
u\|_{L^p((-T,T); L^q(\R^{n}_x; L^2(M^1_y))}\end{equation}
\begin{thm}\label{cauchy}
Let $n\geq 1$ be fixed and $\alpha<4/(n+1)$, then for every initial data
$\varphi\in H^1(\R^n\times M^1)$,
the Cauchy problem \eqref{NLS} has a unique global solution 
$u(t, x, y)$ satisfying : 
$$u(t, x, y) \in {\mathcal C}((-T, T); H^1(\R^n\times M^1)) \cap 
X_{T}\cap Y_{T},
\hbox{ } \forall T>0$$
\end{thm}
\begin{remark}
The main difficulty in the analysis of the Cauchy problem \eqref{NLS}
(compared with the Cauchy problem in the euclidean space)
is related with the fact that
the propagator $e^{-it\Delta_{x,y}}$
on $\R^n\times M^1_y$ does not satisfies the Strichartz estimates which are available
for the propagator
$e^{-it\Delta_{\R^{n+k}}}$ on the euclidean space $\R^{n+k}$.
\end{remark}
Let us now describe some other known cases when \eqref{NLS} is well-posed in $H^1(\R^n\times M^k)$ under the assumption $\alpha<4/(n+k)$.
Using the analysis of \cite{BGT1,BGT2} one may prove such a well-posedness result in the case $\R\times M^2$, i.e. $n=1$ and $k=2$. Moreover, using the analysis of the recent papers \cite{HTT} and \cite{IP} one may also prove such a well-posedness result in the cases $\R^2\times\T^2$ and $\R\times\T^3$ respectively. 
\begin{notation}
Next we fix some notations.
We denote by $L^p_x$ and $H^s_x$ respectively
the space $L^p(\R^{n}_x)$ and $H^s(\R^{n}_x)$.
We also use the notation $L^p_{x,y}=L^p(\R^n_x\times M^k_y)$
and $L^p_xL^q_y=L^p(\R^n_x;L^q(M^k_y))$.
If $v(t)$ is a time dependent function defined on $\R_t$ and valued 
in a Banach space $X$, then we define
$$\|v\|_{L^p_t(X)}^p=\int_\R \|v(t)\|_X^p dt$$
For every $p\in [1, \infty]$ we denote by $p'\in [1, \infty]$ 
its conjugate H\"older exponent.
We denote by $e^{-it\Delta_{x,y}}$ the free propagator associated to the 
Schr\"odinger equation on $\R^{n}_x\times M^k_y$.
\end{notation}

\section{Some useful results on the 
euclidean space $\R^{n}_x$ with $n\geq 1$}\label{cazlions}

In this section we recall some well known facts (see \cite{C})
related to the following minimization problem on $\R^{n}_x$:
\begin{equation}\label{Id-1}
I^\rho_{n,\alpha}=\inf_{\substack{u\in H^1_x\\ \|u\|_{L^2_x}=
\rho}} {\mathcal E}_{n, \alpha} (u)\end{equation}
where for $\alpha<4/n$
\begin{equation}\label{energyd-1}
{\mathcal E}_{n,\alpha}(u)=
\frac 12
\int_{\R^{n}_x}  |\nabla_x u|^2 dx
- \frac 1{2+\alpha}  \int_{\R^{n}_x} |u|^{2+\alpha} dx\end{equation}
By an elementary rescaling argument we have
\begin{equation}\label{potespl}
I^\rho_{n, \alpha}= \rho^{(8+4\alpha-2\alpha n)/(4-\alpha n)} 
I^1_{n,\alpha}\end{equation}
It is well--known that 
\begin{equation}\label{sign}
-\infty<I^\rho_{n, \alpha}<0, \hbox{ } \forall \rho>0
\end{equation}
and
\begin{equation}\label{idfond}
{\mathcal M}^\rho_{n, \alpha}={\mathcal N}_{n,\omega(\rho), \alpha} 
\end{equation}
where ${\mathcal N}_{n,\omega, \alpha}$ is defined in \eqref{Nomega},
\begin{equation}\label{char}
{\mathcal M}_{n,\alpha}^\rho=\{u\in H^1_x|\|u\|_{L^2_x}=\rho \hbox{ and } 
{\mathcal E}_{n, \alpha}(u)=I^\rho_{n,\alpha}\}
\end{equation}
and $\omega(\rho)$ is defined uniquely (see \eqref{remtriv})
by the relation
$$\|u_{n,\omega(\rho), \alpha}\|_{L^2_x}=\rho$$
We also recall that the functions $u_{n,\omega, \alpha}$ (defined as the unique radially symmetric
and positive solution to \eqref{uomegaintro}) 
satisfy the following Pohozaev type identity
(for a proof of \eqref{poza} see the proof of \eqref{poho} in next section):
\begin{equation}\label{poza}
\int_{\R^{n}_x} |\nabla_x u_{n,\omega, \alpha}|^2 dx= 
\frac{\alpha n}{2(\alpha+2)} \int_{\R^{n}_x}
|u_{n,\omega, \alpha}|^{2+\alpha} dx\end{equation}
On the other hand if we multiply \eqref{uomegaintro} by $u_{n,\omega,\alpha}$ and we integrate by parts
then we get
$$\int_{\R^{n}_x} |\nabla_x u_{n,\omega, \alpha}|^2 dx
+\omega \|u_{n,\omega,\alpha}\|_{L^2_x}^2= 
\int_{\R^{n}_x}
|u_{n,\omega, \alpha}|^{2+\alpha} dx$$
that in conjunction with \eqref{poza}
gives
\begin{equation}\label{pozae}
\omega \|u_{n,\omega, \alpha}\|_2^2= \frac {2\alpha + 4 - \alpha n}{\alpha n} \int_{\R^{n}_x} 
|\nabla_x u_{n,\omega, \alpha}|^2 dx
\end{equation}
$$=\frac {4\alpha + 8 - 2\alpha n}{\alpha n-4} \Big ( \frac 12 \int_{\R^{n}_x} |\nabla_x u_{n,\omega, \alpha}|^2 dx 
-\frac{1}{2+ \alpha}
\int_{\R^{n}_x} |u_{n,\omega, \alpha}|^{2+\alpha} dx\Big )
$$$$=\frac {4\alpha + 8 - 2\alpha n}{\alpha n-4} I_{n, \alpha}^{\|u_{n,\omega, \alpha}\|_{L^2_x}}$$
(at the last step we have used the fact that due to \eqref{idfond}
we have that $u_{n,\omega, \alpha}$ is a minimizer for 
${\mathcal E}_{n, \alpha}$ on its associated constrained).\\
Finally notice that by \eqref{poza}
we deduce 
\begin{equation}\label{d-1min}
I_{n, \alpha}^{\|u_{n,\omega, \alpha}\|_{L^2_x}}={\mathcal E}_{n, \alpha}(u_{n,\omega, \alpha})=
\frac{\alpha n -4}{2\alpha n} \int_{\R^{n}_x} 
|\nabla_x u_{n,\omega, \alpha}|^2 dx
\end{equation}

\section{An auxiliary problem}

In this section we study the minimizers of the following
minimization problems
\begin{equation}\label{JnMk}J_{n, M^k, \alpha,\lambda}= \inf_{\substack{u\in H^1(\R^n \times M^k)\\ \|u\|_{L^2_{x,y}}=1}} 
{\mathcal E}_{n, M^k, \alpha, \lambda} (u)\end{equation} where
$$ {\mathcal E}_{n, M^k, \alpha, \lambda}(u)=
\int_{M^k_y}\int_{\R^{n}_x}  \Big (\frac \lambda 2 |\nabla_y u|^2 + \frac 12 |\nabla_x u|^2
- \frac 1{2+\alpha} |u|^{2+\alpha}\Big ) dxdvol_{M_y^k}$$ 
We also introduce the following sets:
$${\mathcal M}_{n, M^k,\alpha, \lambda}
=\{w\in H^1(\R^n \times M^k) |\|w\|_{L^2_{x,y}}=1 \hbox{ and } \mathcal{E}_{n, M^k, \alpha, \lambda}(w)=
J_{n, M^k, \alpha, \lambda}\}$$
\begin{thm}\label{rigidity}
Let $n, M^k$ and $0<\alpha<\frac 4{n+k}$ be given. There exists $\lambda^*\in (0, \infty)$
such that: 
\begin{equation}\label{prop12}{\mathcal M}_{n, M^k, \alpha, \lambda}
={\mathcal N}_{n,\bar \omega, \alpha}, \hbox{ } \forall \lambda>\lambda^*\end{equation}
and
\begin{equation}\label{prop14}{\mathcal M}_{n, M^k, \alpha, \lambda}\cap
{\mathcal N}_{n,\bar \omega, \alpha}=\emptyset, \hbox{ } \forall \lambda<\lambda^*
\end{equation}
where
$\bar \omega$ is defined by the condition
$$vol(M^k) \|u_{n,\bar \omega, \alpha}\|_{L^2_x}^2=1$$
\end{thm}
We 
fix a sequence $\lambda_j\rightarrow \infty$ and a corresponding sequence of functions 
$u_{\lambda_j} \in {\mathcal M}_{n, M^k,\alpha,\lambda_j}$.
In the sequel we shall 
assume that
\begin{equation}\label{positivity}
u_{\lambda_j}(x,y)\geq 0, \hbox{ } 
\forall (x,y)\in \R_x^{n}\times M^k_y
\end{equation}
Indeed it is well-known that if $u_{\lambda_j}$ is a minimizer, then also
$|u_{\lambda_j}|$ is a minimizer. In particular there exists
at least one minimizer
which satisfies \eqref{positivity}.
\\
Notice that the functions $u_{\lambda_j}$ depend
in principle on the full set of variables $(x, y)$. Our aim is to prove
that for $j$ large and up to subsequence,
the functions $u_{\lambda_j}$ will not depend explicitly
on the variable $y$.\\
First we prove some a priori bounds
satisfied by $u_{\lambda_j}(x, y)$.
Recall that the quantities $I^\rho_{n,\alpha}$ are defined in \eqref{Id-1}.
\begin{lem} Assume the same assumptions as in Theorem \ref{rigidity},
then we have:
\begin{equation}\label{limite}\lim_{j \rightarrow \infty} 
J_{n, M^k, \alpha,\lambda_j}= vol (M^k)
I^{1/\sqrt {vol(M^k)}}_{n, \alpha}
\end{equation}
and
\begin{equation}\label{improved}
\lim_{j\rightarrow \infty} \lambda_j \int_{M^k_y} \int_{\R^{n}_x}
|\nabla_y u_{\lambda_j} |^2 dxdvol_{M_y^k}=0
\end{equation}
\end{lem}
{\bf Proof.}
First notice that
\begin{equation}\label{minim}
J_{n, M^k, \alpha,\lambda_j}\leq vol (M^k)
I^{1/\sqrt {vol(M^k)}}_{n, \alpha}
\end{equation}
In fact let $w(x) \in H^1_x$ be such that $\|w\|_{L^2_x}=\frac{1}{\sqrt {vol(M^k)}}
$ and
${\mathcal E}_{n, \alpha}(w)=I^{1/\sqrt {vol(M^k)}}_{n, \alpha}$.
Then
we get easily:
$$J_{n, M^k, \alpha,\lambda_j} \leq {\mathcal E}_{n, M^k, \alpha,\lambda_j} (w(x))$$$$
=vol (M^k)  \Big(\frac 12
\int_{\R^{n}_x}  |\nabla_x w|^2 dx - 
\frac 1{2+ \alpha}  \int_{\R^{n}_x} |w|^{2+\alpha} dx\Big )
= vol(M^k)  I^{1/\sqrt {vol(M^k)}}_{n,\alpha}$$
and this concluded the proof of \eqref{minim}.\\
Next we claim that
\begin{equation}\label{partial}
\lim_{j\rightarrow \infty} \int_{M^k_y} \int_{\R^{n}_x}
|\nabla_y u_{\lambda_j}|^2 dxdvol_{M_y^k}=0\end{equation}
In order to prove this fact assume by the absurd that it is false
then
there exists a subsequence of $\lambda_j$ (that we 
still denote by $\lambda_j$)
such that
$$\lim_{j\rightarrow \infty} \lambda_j=\infty
\hbox{ and } 
\int_{M^k_y} \int_{\R_x^{n}}|\nabla_y u_{\lambda_j}|^2 dxdvol_{M^k_y}
\geq \epsilon_0>0$$
and in particular
\begin{equation}\label{infili}
\lim_{j\rightarrow \infty} (\lambda_j -1)\int_{M^k_y} \int_{\R_x^{n}}
|\nabla_y u_{\lambda_j}|^2 dxdvol_{M^k_y}
=\infty\end{equation}
On the other hand by the classical Gagliardo Nirenberg inequality (see \eqref{GagN})
we deduce the existence of $0<\mu<2$ such that: 
$$\frac 12 \int_{M^k_y} \int_{\R_x^{n}} 
(|\nabla_y v|^2 + |\nabla_x v|^2
+ |v|^2) dx dvol_{M_y^k}   
- \frac{1}{2+\alpha} \int_{M^k_y} \int_{\R_x^{n}} |v|^{2+\alpha}dxdvol_{M_y^k} 
$$$$\geq 
\frac 12 \int_{M^k_y} \int_{\R_x^{n}} 
(|\nabla_y v|^2 + |\nabla_x v|^2
+ |v|^2) dx dvol_{M^k_y}  $$$$- C 
\left [\int_{M^k_y} \int_{\R_x^{n}} 
(|\nabla_y v|^2 + |\nabla_x v|^2
+ |v|^2) dx dvol_{M_y^k}\right]^\mu$$$$\geq \inf_{t>0} (1/2 t^2 - C t^\mu)=C(\mu)>-\infty  
$$ 
$$\forall v\in H^1(\R^n \times M^k) \hbox{ s.t. } \|v\|_{L^2_{x,y}}=1$$
By the previous inequality we get
$$\mathcal E_{n, M^k, \alpha, \lambda_j} (v) - 
\frac 12 (\lambda_j-1) \int_{M^k_y} \int_{\R^{n}_x}
|\nabla_y v|^2 \geq - \frac 12 + C(\mu) \hbox{ }$$$$
\forall v\in H^1(\R^n \times M^k) \hbox{ s.t. } \|v\|_{L^2_{x,y}}=1$$
In particular if we choose $v=u_{\lambda_j}$
then we get
$$J_{n, M^k, \alpha, \lambda_j}= \mathcal E_{n, M^k, \alpha, \lambda_j}(u_{\lambda_j})
$$$$\geq \frac 12 (\lambda_j-1) \int_{M^k_y} \int_{\R^{n}_x} 
|\nabla_y u_{\lambda_j}|^2 dxdvol_{M_y^k} - 
\frac 12 + C(\mu)$$
By \eqref{infili} this implies $\lim_{n\rightarrow \infty}
J_{n, M^k, \alpha,\lambda_j}=\infty$ 
and this is in contradiction with \eqref{minim}. Hence \eqref{partial}
is proved.\\
\\
Next we introduce the functions
$$w_j(y)=\|u_{\lambda_j}(x, y)\|_{L^2_x}^2$$
Notice that 
\begin{equation}\label{L1}\|w_j(y)\|_{L^1_y}=1
\end{equation}
and moreover
$$\int_{M^k_y} |\nabla_y w_j(y)|dvol_{M^k_y}
\leq C \int_{M^k_y}\int_{\R^n_x} |u_{\lambda_j}(x, y)| 
|\nabla_y u_{\lambda_j}(x, y)| dxdvol_{M^k_y}$$
$$\leq C \|u_{\lambda_j}\|_{L^2_{x,y}}\|\nabla_y u_{\lambda_j}\|_{L^2_{x,y}}$$
Hence due to \eqref{partial} we get 
\begin{equation}\label{L1grad}\lim_{j\rightarrow \infty}\|\nabla_y 
w_j\|_{L^1_y}=0\end{equation}
By combining \eqref{L1} and \eqref{L1grad} with the Rellich compactness theorem 
and with the Sobolev embedding 
$W^{1,1}(M^1)\subset L^\infty(M^1)$ and $W^{1,1}(M^2)\subset L^2(M^2)$
we deduce respectively in the case $k=1$ and $k=2$
that (up to a subsequence)
\begin{equation}\label{unifcont1}\lim_{j\rightarrow \infty}
\|w_j(y)-1/vol(M^1)\|_{L^r_y}=0,
\hbox{ } \forall 1\leq r <\infty
\end{equation}
and
\begin{equation}\label{unifcont2}\lim_{j\rightarrow \infty}
\|w_j(y)-1/vol(M^2)\|_{L^r_y}=0,
\hbox{ } \forall 1\leq r <2
\end{equation}
For $k>2$ we use the Sobolev embedding $H^1(M^k)\subset L^{2k/(k-2)}(M^k)$
and we get
$$\sup_j \|u_{\lambda_j}\|_{L^2_x L^{2k/(k-2)}_y}\leq C \sup_j 
\|u_{\lambda_j}\|_{L^2_x H^1(M^k_y)}<\infty$$
(where at the last step we have used the  fact $\sup_j \big(\|u_{\lambda_j}\|_{L^2_{x,y}}
+\|\nabla_y u_{\lambda_j}\|_{L^2_{x,y}}\big)<\infty$).
By the Minkowski inequality the bound above implies
$\sup_j \|u_{\lambda_j}\|_{L^{2k/(k-2)}_yL^2_x}$ which is equivalent to the condition
\begin{equation}\label{maxSob}\sup_j  \|w_j(y)\|_{L^{k/(k-2)}_y}<\infty \hbox{ for } k>2\end{equation}
By combining \eqref{L1} and \eqref{L1grad} with the Rellich compactness theorem 
we deduce that up to a subsequence
$$\|w_j(y)-1/vol(M^k)\|_{L^1_y}=0 \hbox{ for } k>2$$
and hence by interpolation with \eqref{maxSob} we get
\begin{equation}\label{unifcontk}
\|w_j(y)-1/vol(M^k)\|_{L^r_y}=0 \hbox{ for } k>2, 1\leq r<k/(k-2)
\end{equation}
By the definition of 
$I^\rho_{n,\alpha}$ (see \eqref{Id-1}) and \eqref{potespl} we get
\begin{equation}\label{fixedy}\frac 12 \int_{\R^{n}_x} |\nabla_x u_{\lambda_j}(x,y)|^2 dx
- \frac 1{2+\alpha}\int_{\R^{n}_x} |u_{\lambda_j} (x,y)
|^{2+\alpha} dx\end{equation}
$$\geq I_{n, \alpha}^{\|u_{\lambda_j}(\cdot,y)\|_{L^2_x}}=I^1_{n, \alpha} \|u_{\lambda_j}(\cdot,y)
\|_{L^2_x}^{(8+4\alpha-2\alpha n)/(4-\alpha n)}=I^1_{n, \alpha} w_j(y)^{(4+2\alpha-\alpha n)/(4-\alpha n)}
$$
$$\hbox{ } \forall y\in M^k,
\hbox{ } \forall j\in\N$$
Next notice that by definition
\begin{equation}\label{massi}
J_{n, M^k, \alpha,\lambda_j}={\mathcal E}_{n, M^k, \alpha,\lambda_j}(u_{\lambda_j})\end{equation}
$$= \frac 12 \int_{M^k_y} \int_{\R_x^{n}} (\lambda_j 
|\nabla_y u_{\lambda_j}|^2 + 
|\nabla_x u_{\lambda_j}|^2 )dxdy
- \frac{1}{2+\alpha} \int_{M^k_y} \int_{\R_x^{n}} |u|^{2+\alpha} dxdvol_{M^k_y}
$$
and we can continue 
\begin{equation}\label{massi2}...\geq \int_{M^k_y} \Big (   
\frac 12 \int_{\R^{n}_x} |\nabla_x u_{\lambda_j}( x,y)|^2 dx
- \frac 1{2+\alpha}\int_{\R^{n}_x} |u_{\lambda_j} (x,y)
|^{2+\alpha} dx\Big ) dvol_{M^k_y}\end{equation}
$$\geq I^1_{n,\alpha} \int_{M^k_y} w_j(y)^{(4+2\alpha-\alpha n)/(4-\alpha n)}dvol_{M_y^k}
$$$$= I^1_{n,\alpha} vol(M^k) vol(M^k)^{-(4+2\alpha-\alpha n)/(4-\alpha n)}+o(1)$$
where $o(1)\rightarrow 0$ as $j\rightarrow \infty$ and at the last step we have combined
\eqref{unifcont1}, \eqref{unifcont2} and \eqref{unifcontk} respectively for $k=1$, $k=2$ 
and $k>2$ and we used our assumption on $\alpha$.
By combining this fact with 
\eqref{potespl} we have 
\begin{equation}\label{matt}\liminf_{j\rightarrow \infty }
J_{n, M^k, \alpha,\lambda_j}\geq 
vol(M^k)  I^{1/\sqrt {vol(M^k)}}_{n, \alpha}\end{equation}
Hence
\eqref{limite} follows by combining \eqref{minim}
with \eqref{matt}.\\
\\
Next we prove \eqref{improved}. For that purpose, it suffices to keep the term $\lambda_j 
|\nabla_y u_{\lambda_j}|^2$ in the previous analysis. Namely, by combining \eqref{limite}
with \eqref{massi} and \eqref{massi2}
we get
\begin{equation}\label{utile}
vol(M^k) I_{n, \alpha}^{1/\sqrt {vol(M^k)}}+ 
g(j)\geq \frac 12 \lambda_j \int_{M^k_y} \int_{\R^{n}_x} 
|\nabla_y u_{\lambda_j}|^2 dx dvol_{M^k_y} + h(j)\end{equation}
where
$$\lim_{j\rightarrow \infty}g(j)= 0$$
and
$$\liminf_{j\rightarrow \infty} h(j)\geq vol(M^k) 
I_{n, \alpha}^{1/\sqrt {vol(M^k)}}$$
Hence \eqref{improved} follows by \eqref{utile}.

\hfill$\Box$

\begin{lem}
We have the following identity:
\begin{equation}\label{poho}
\int_{M^k_y} \int_{\R_x^{n}} |\nabla_x u_{\lambda_j}|^2 dxdvol_{M^k_y}  
=\frac{\alpha n}{2(2+\alpha)}
\int_{M^k_y}\int_{\R_x^{n}}|u_{\lambda_j}|^{2+\alpha} dxdvol_{M^k_y}
\end{equation}
Moreover there exist $J\in \N$ such that
$$\forall j>J \hbox{ }\exists \omega(\lambda_j)>0 \hbox{ s.t. }$$
\begin{equation}\label{lagra}
-\lambda_j \Delta_y u_{\lambda_j} - \Delta_x u_{\lambda_j} +
\omega(\lambda_j) u_{\lambda_j}= u_{\lambda_j} |u_{\lambda_j}|^{\alpha}
\end{equation}
and the following limit exists
\begin{equation}\label{Lim}
\lim_{j \rightarrow \infty}
\omega(\lambda_j)=\bar \omega\in (0, \infty)
\end{equation}
\end{lem}

{\bf Proof.} Since $u_{\lambda_j}$ is a constrained minimizer
for ${\mathcal E}_{n,M^k,\alpha,\lambda_j}$ on the ball of size
$1$ in $L^2(\R^n\times M^k)$, then we get 
$$\frac d{d\epsilon} \Big [
{\mathcal E}_{n, M^k, \alpha,\lambda_j} (\epsilon^{\frac{n}{2}} 
u_{\lambda_j} (\epsilon x,y)
\Big ]_{\epsilon=1}=0$$ 
which is equivalent to 
$$\frac d{d\epsilon}  \Big [\frac 12 \lambda_j \int_{M^k_y} \int_{\R^{n}_x} 
|\nabla_y u_{\lambda_j}|^2 dxdvol_{M^k_y} $$
$$ + \frac 12 \epsilon^2 \int_{M^k_y} \int_{\R^{n}_x} 
|\nabla_x u_{\lambda_j}|^2 dxdvol_{M^k_y} 
- \frac 1{2+\alpha} \epsilon^{\alpha n/2}\|u_{\lambda_j}
\|_{L^{2+\alpha}_{x,y}}^{2+\alpha}\Big ]_{\epsilon=1}=0$$
By computing explicitly the derivative 
(in $\epsilon$) we deduce \eqref{poho}.\\
\\
Next notice that by using the Lagrange multiplier technique
we get \eqref{lagra} for a suitable $\omega(\lambda_j)\in \R$.
On the other hand by \eqref{lagra}
we get
$$\int_{M^k_y} \int_{\R^{n}_x} 
(\lambda_j |\nabla_y u_{\lambda_j}|^2 +|\nabla_x u_{\lambda_j}|^2) dxdvol_{M^k_y}
+ \omega(\lambda_j) \|u_{\lambda_j}\|_{L^2_{x,y}}^2
$$$$= \int_{M^k_y} \int_{\R^{n}_x} |u_{\lambda_j}|^{2+\alpha} dxdvol_{M_y^k}$$
that by \eqref{poho} gives 
$$
\omega(\lambda_j) = \frac {-\alpha n + 4+2\alpha}{\alpha n} 
\int_{M^k_y} \int_{\R^{n}_x}  |\nabla_x u_{\lambda_j}|^2 dxdvol_{M_y^k} 
-\lambda_j \int_{M_y^k} \int_{\R^{n}_x}  |\nabla_y u_{\lambda_j}|^2 dxdvol_{M_y^k}$$
and hence by \eqref{improved} we get
\begin{equation}\label{moltiplicatore}
\omega(\lambda_j) = \frac {-\alpha n + 4+2\alpha}{\alpha n}
\int_{M^k_y} \int_{\R^{n}_x}  |\nabla_x u_{\lambda_j}|^2 dxdvol_{M_y^k} 
+ o(1)\end{equation}
where $\lim_{j\rightarrow \infty} o(1)=0$.\\
On the other hand notice that by \eqref{poho}
we get
$$J_{n, M^k, \alpha,\lambda_j}= {\mathcal E}_{n, M^k, \alpha,\lambda_j}(u_{\lambda_j})$$$$=
\frac {\alpha n-4}{2\alpha n}  \int_{M^k_y} \int_{\R^{n}_x}  |\nabla_x u_{\lambda_j}|^2 dxdvol_{M_y^k}
+ \frac 12 \int_{M^k_y} \int_{\R^{n}_x} \lambda_j |\nabla_y u_{\lambda_j}|^2 dxdvol_{M_y^k}$$
and 
by \eqref{improved} 
\begin{equation}\label{lastnught}\int_{M_y^k} \int_{\R^{n}_x}  
|\nabla_x u_{\lambda_j}|^2 dxdvol_{M_y^k}= 
\frac{2\alpha n}{\alpha n -4} J_{n, M^k, \alpha,\lambda_j}+o(1)\end{equation}
By \eqref{limite} it implies
\begin{equation}\label{infimum}
\int_{M^k_y} \int_{\R^{n}_x}  |\nabla_x u_{\lambda_j}|^2 dxdvol_{M_y^k}= 
\frac{2\alpha n}{\alpha n-4} vol(M^k) 
I_{n,\alpha}^{1/\sqrt {vol(M^k)}}+ o(1)\end{equation}
that in conjunction with \eqref{moltiplicatore} and \eqref{sign}
implies $\omega(\lambda_j)>0$  for $j$ large
enough. Moreover \eqref{Lim} follows by \eqref{moltiplicatore} 
and \eqref{infimum}.

\hfill$\Box$

Next recall that the sets ${\mathcal M}^\rho_{n, \alpha}$ are the ones defined in \eqref{char}.
\begin{lem}\label{firstrig}
Let $\bar \omega$ be as in \eqref{Lim} and let
$v(x)\in {\mathcal M}_{n,\alpha}^{1/\sqrt {vol(M^k)}}$ be such that
$v(x)>0$.
Then 
$$-\Delta_x v + \bar \omega v= v|v|^{\alpha}$$ 
\end{lem}

{\bf Proof.}
It is well-known that $$-\Delta_x v + \omega_1 v= v|v|^{\alpha}$$
for a suitable $\omega_1>0$. More precisely we 
can assume that up to translation
$v=u_{n,\omega_1, \alpha}$. 
Our aim is to prove that $\omega_1=\bar \omega$. 
Notice that by \eqref{pozae}
\begin{equation}\label{omega1}
\omega_1 \frac{1}{vol(M^k) }= \frac {4\alpha + 8 - 2\alpha n}{\alpha n-4} 
I_{n, \alpha}^{\|v\|_{L^2_x}}=
\frac {4\alpha + 8 - 2\alpha n}{\alpha n-4}I^{1/\sqrt {vol(M^k)}}_{n, \alpha}
\end{equation}
On the other hand by 
\eqref{moltiplicatore} and \eqref{infimum} we get
$$\omega(\lambda_j) = \frac {-2\alpha n + 8+4\alpha}{\alpha n-4}vol(M^k) 
I_{n, \alpha}^{1/\sqrt {vol(M^k)}}+ o(1)$$
and hence passing to the limit in $j$ we get
\begin{equation}\label{omegabar}
\bar \omega = \frac {-2\alpha n + 8+4\alpha}{\alpha n-4} 
vol(M^k) I_{n, \alpha}^{1/\sqrt {vol(M^k)}}
\end{equation}
By combining \eqref{omega1} and \eqref{omegabar}
we get $\bar \omega=\omega_1$.

\hfill$\Box$

\begin{lem}\label{rescatr}
There exist a subsequence of $\lambda_j$ (that we shall denote still by $\lambda_{j}$)
and a sequence $\tau_j\in \R^{n}_x$
such that
$$\lim_{j \rightarrow \infty} \|u_{\lambda_{j}} 
(x+ \tau_j, y)-u_{\bar \omega}\|_{H^1(\R^n\times M^k)}=0$$
where $u_{\bar \omega}\in {\mathcal N}_{n,\bar \omega, \alpha}$,
$u_{\bar \omega}>0$ and $\bar \omega$
is defined in \eqref{Lim}.
\end{lem}

{\bf Proof.} 
By combining \eqref{improved} and \eqref{infimum},
and  since $\|u_{\lambda_j}\|_{L^2_{x,y}}=1$, we deduce that
$u_{\lambda_j}$ is bounded in $H^1(\R^n\times M^k)$.
Moreover by combining \eqref{limite}
with the fact that $I_{n,\alpha}^{1/\sqrt {vol(M^k)}}<0$ (see \eqref{sign})
then we get 
$$\inf_{j} \|u_{\lambda_j}\|_{L^{2+\alpha}_{x,y}}>0$$
By using the localized version of
the Gagliardo Nirenberg inequality \eqref{GNprecised} 
(in the same spirit as in the appendix) we get
the existence (up to subsequence) of $\tau_j\in \R^{n}_x$ such that 
$$u_{\lambda_{j}}(x+\tau_j, y) \rightharpoonup w\neq 0 \hbox{ in } H^1(\R^n\times M^k)$$
Moreover
due to \eqref{positivity}
we can assume that 
$$w(x, y)\geq 0 \hbox{ a.e. } (x, y)\in \R^{n}_x\times M^k_y$$
and by \eqref{improved}
we get $\nabla_y w=0$.
In particular $w$ is $y$-independent.\\
By combining
\eqref{improved} and \eqref{Lim} we pass 
to the limit in \eqref{lagra} 
in the distribution sense and we get
\begin{equation}\label{eul-lagr}
-\Delta_x w +\bar \omega w=w|w|^{\alpha} \hbox{ in } \R_x^{n}, \hbox{ } w(x)\geq 0, 
\hbox{ } w\neq 0
\end{equation}
We claim that 
\begin{equation}\label{strongL2}\|w\|_{L^2_x}=\frac{1}{\sqrt {vol(M^k)}}
\end{equation}
If not then we can assume
$\|w\|_{L^2_x}=\beta<\frac{1}{\sqrt {vol(M^k)}}$
and since $w$ solves \eqref{eul-lagr} by
\eqref{idfond}
we get
\begin{equation}\label{absalpha}
w\in {\mathcal M}_{n,\alpha}^\beta
\end{equation}
On the other hand 
by Lemma \ref{firstrig} the equation \eqref{eul-lagr}
is satisfied by any $v\in {\mathcal M}_{n, \alpha}^\frac{1}{\sqrt {vol(M^k)}}$.
Hence again by \eqref{idfond} and by the injectivity
of the map $\rho\rightarrow \omega(\rho)$
(see \eqref{remtriv}) we deduce that necessarily
$\beta= \frac{1}{\sqrt {vol(M^k)}}$.\\
In particular by \eqref{strongL2} we deduce
$$\lim_{j\rightarrow \infty} \|u_{\lambda_{j}}(x+\tau_j,y) - w\|_{L^2_{x,y}}=0$$
Next notice that
by \eqref{improved} and since we have already proved that $\nabla_y w=0$ we deduce that
$$\lim_{j\rightarrow \infty} \|\nabla_y u_{\lambda_{j}}(x+\tau_j,y)\|_{L^2_{x,y}}
=0=\|\nabla_y w\|_{L^2_{x,y}}$$
Hence in order to conclude that
$u_{\lambda_{j}}(x+\tau_j,y)$ converges strongly to $w$ in $H^1(\R^n \times M^k)$
it is sufficient to prove
that
$$\lim_{j\rightarrow \infty} \|\nabla_x u_{\lambda_{j}}(x+\tau_j,y)\|_{L^2_{x,y}}
=\sqrt {vol(M^k)}\|\nabla_x w\|_{L^2_x}=\|\nabla_x w\|_{L^2_{x,y}}$$
This last fact follows by combining \eqref{d-1min} (where we use the fact that $w\in {\mathcal N}_{n,\bar \omega,
\alpha}$
by \eqref{eul-lagr}
and $\|w\|_{L^2_x}=\frac{1}{\sqrt {vol(M^k)}}$ by \eqref{strongL2})
and \eqref{infimum}.

\hfill$\Box$

\begin{lem}\label{unod}
There exists $j_0>0$ such that
$$\nabla_y u_{\lambda_j} =0, \hbox{ } \forall j>j_0$$
\end{lem}

{\bf Proof.}
By Lemma \ref{rescatr} we can assume that
\begin{equation}\label{convstrong}
u_{\lambda_j}\rightarrow u_{\bar \omega}
\hbox{ in } H^1(\R^n\times M^k)\end{equation}
We introduce $w_j=\sqrt{-\Delta_y} u_{\lambda_j}$. Notice that due to \eqref{lagra}
the functions $w_j$ satisfy 
\begin{equation}\label{minc}-\lambda_j \Delta_y w_j - \Delta_x w_j +
\omega(\lambda_j) w_j= \sqrt{-\Delta_y} (u_{\lambda_j} |u_{\lambda_j}|^{\alpha})
\end{equation}
that after multiplication by $w_j$ implies 
\begin{equation}\label{contrad}
\int_{M^k_y} \int_{\R^{n}_x} 
\Big[
\lambda_j |\nabla_y w_j |^2 
+|\nabla_x w_j|^2 + \omega(\lambda_j) |w_j|^2 
\end{equation}
$$- \sqrt{-\Delta_y} (u_{\lambda_j} |u_{\lambda_j}|^{\alpha})w_j\Big] dxdvol_{M_y^k}=0
$$
In turn it gives
\begin{equation}\label{zerobas}0=\int_{M^k_y} \int_{\R^{n}_x} (\lambda_j-1)
|\nabla_y w_j|^2
- (\alpha+1) \sqrt{-\Delta_y} (u_{\lambda_j} |u_{\bar \omega}|^{\alpha})w_j
dxdvol_{M_y^k}+\end{equation}
$$ \int_{M^k_y} \int_{\R^{n}_x} (|\nabla_y w_j|^2 + | \nabla_x w_j|^2 
+ \bar \omega   |w_j|^2 +
\sqrt{-\Delta_y} (u_{\lambda_j} ((\alpha+1)|u_{\bar \omega}|^\alpha - 
|u_{\lambda_j}|^{\alpha})) w_j dxdvol_{M_y^k}
$$
$$+\int_{M^k_y} \int_{\R^{n}_x} (\omega(\lambda_j) - \bar \omega) |w_j|^2 dxdy\equiv$$
$$I_j+II_j+III_j$$
Next we fix an orthonormal basis of eigenfunctions for $-\Delta_y$, i.e.
$-\Delta_y \varphi_k =\mu_k\varphi_k$ and $\varphi_0={\rm const}$.
We can write the following development
\begin{equation}\label{fouriel}
w_j(x, y)=\sum_{k\in \N\setminus \{0\}} a_{j, k}(x)\varphi_k(y)
\end{equation}
(where the eigenfunction $\varphi_0$ does not enter in the development).
By using the representation 
in \eqref{fouriel} 
we get
\begin{equation}\label{Ipet}I_j\geq \sum_{k\neq 0} (\lambda_j -1)|\mu_k|^2
\int_{\R^{n}_x} |a_{j,k}(x)|^2 dx - 
(\alpha+1)\sum_{k\neq 0} \int_{\R^{n}_x} |u_{\bar \omega}(x)|^{\alpha} |a_{j,
k}(x)|^2 dx
\end{equation}
and by \eqref{Lim} we get 
\begin{equation}\label{III} III_j=o(1)\|w_j\|_{L^2_{x,y}}^2\end{equation}
By combining \eqref{Ipet} with \eqref{III} we get
\begin{equation}\label{I+III}
I_j+III_j\geq 0\end{equation}
for $j$ large enough.
In order to estimate $II_j$ notice that
by the Cauchy-Schwartz inequality
we get 
\begin{equation}\label{strof}
\Big|\int_{M_y^k}\int_{\R^n_x} \sqrt{-\Delta_y} (u_{\lambda_j} ((\alpha+1)|u_{\bar
\omega}|^\alpha - 
|u_{\lambda_j}|^{\alpha}))w_j dxdvol_{M_y^k}\Big|
\end{equation}
$$\leq \|\sqrt{-\Delta_y} (u_{\lambda_j} ((\alpha+1)|u_{\bar \omega}|^\alpha - 
|u_{\lambda_j}|^{\alpha}))\|_{L^{\frac{2(n+k)}{n+k+2}}_x L^{\frac{2(n+k)}{n+k+2}}_y}
\|w_j\|_{L^{\frac{2(n+k)}{n+k-2}}_{x,y}}
$$
$$\leq C \|\nabla_y (u_{\lambda_j} ((\alpha+1)|u_{\bar \omega}|^\alpha - 
|u_{\lambda_j}|^{\alpha}))\|_{L^{\frac{2(n+k)}{n+k+2}}_x L^{\frac{2(n+k)}{n+k+2}}_y}
\|w_j\|_{L^{\frac{2(n+k)}{n+k-2}}_{x,y}}
$$
where at the last step we have used the following estimate
\begin{equation}\label{equiv}
\forall p\in (1, \infty) \hbox{ } \exists c(p), C(p)>0 \hbox{ s.t. }\end{equation}
$$c(p) \|\sqrt{-\Delta_y} f\|_{L^p_y}
\leq \|\nabla_y f\|_{L^p_y}\leq C(p) \|\sqrt{-\Delta_y} f\|_{L^p_y}
$$
Indeed, using \cite[Theorem 3.3.1]{Sogge}, we have that $\sqrt{-\Delta_y}$ is a first order classical pseudo differential operator on $M$ with a principle symbol 
$(g^{i,j}(y)\xi_i\,\xi_j)^{1/2}$.
Observe that 
$$
C_1\sum_{i,j}g^{i,j}(y)\xi_i\,\xi_j\leq \sum_{i} |\sum_{j}g^{i,j}(y)\xi_j|^2 
\leq C_2|\xi|^2\leq C_3 \sum_{i,j}g^{i,j}(y)\xi_i\,\xi_j
$$
Moreover one can assume that in  \eqref{equiv} $f$ has no zero frequency. 
Then one can deduce \eqref{equiv} by working in local coordinates, introducing a classical angular partition of unity according to the index $l\in [1,\cdots,k]$ such that 
$$
\sum_{i,j}g^{i,j}(y)\xi_i\,\xi_j\leq c 
|\sum_{j}g^{l,j}(y)\xi_j|^2 
$$ 
and, most importantly, using the $L^p$ boundedness of zero order 
pseudo differential operators on $\R^k$ (for the proof of this fact we refer to \cite[Theorem 3.1.6]{Sogge}).

Next, by the chain rule we get
$$\nabla_y \Big (u_{\lambda_j} ((\alpha+1)|u_{\bar \omega}|^\alpha - 
|u_{\lambda_j}|^{\alpha})\Big)$$$$= (\alpha +1)
\nabla_y  u_{\lambda_j} \Big (|u_{\bar \omega}|^\alpha - 
|u_{\lambda_j}|^{\alpha}\Big )$$
and by the H\"older inequality we can continue the estimate
\eqref{strof} as follows
$$...\leq C \Big \| \|\nabla_y u_{\lambda_j}\|_{L^q_y} \||u_{\bar \omega}|^\alpha - 
|u_{\lambda_j}|^{\alpha}\|_{L^r_y}\Big \|_{L^{\frac {2(n+k)}{n+k+2}}_x} \|w_j\|_{L^{\frac{2(n+k)}{n+k-2}}_{x,y}}$$
where
$$\frac 1q+\frac 1r=\frac {n+k+2}{2(n+k)}$$
and again by the H\"older inequality in the $x$-variable
we can continue
$$...\leq C \|\nabla_y u_{\lambda_j}\|_{L^q_{x,y}} \||u_{\bar \omega}|^\alpha - 
|u_{\lambda_j}|^{\alpha}\|_{L^r_{x,y}} \|w_j\|_{L^{\frac{2(n+k)}{n+k-2}}_{x,y}}$$
Notice that if we fix 
$$q=\frac{2(n+k)}{n+k-2} \hbox{ and } r=\frac{n+k}{2}$$
then by combining the Sobolev embedding
\begin{equation}\label{crit}H^1_{x,y}\subset L^{\frac{2(n+k)}{n+k-2}}_{x,y}
\end{equation}
with \eqref{convstrong} and \eqref{equiv}, we can continue the estimate
$$
...\leq o(1) \|\sqrt{-\Delta_y} 
u_{\lambda_j}\|_{L^q_{x,y}} \|w_j\|_{H^1_{x,y}}= o(1) \|w_j\|_{H^1_{x,y}}^2 
$$
where $\lim_{j\rightarrow \infty} o(1)=0$.
By combining this information in conjunction with the structure of $II_j$
we get
\begin{equation}\label{IIj}
II_j\geq \|w_j\|_{H^1_{x,y}}^2(1-o(1))\geq 0 \hbox{ for } j>j_0
\end{equation}
By combining \eqref{zerobas}, \eqref{I+III} and \eqref{IIj} we deduce $w_j=0$ for
$j$ large enough.

\hfill$\Box$
\\
\\

{\bf Proof of Theorem \ref{rigidity}}
By using the diamagnetic inequality we deduce that (up to a remodulation factor
$e^{i\theta}$) we can assume that 
$v\in {\mathcal M}_{n,M^k, \alpha, \lambda}$ is real valued.
Moreover if $v\in  {\mathcal M}_{n,M^k, \alpha, \lambda}$ 
then also $|v|\in {\mathcal M}_{n,M^k, \alpha, \lambda}$.
By a standard application of the strong maximum principle we finally 
deduce that it is not restrictive to assume that 
$v\in  {\mathcal M}_{n,M^k, \alpha, \lambda}$ and $v(x, y)>0, \hbox{ }
\forall (x, y)\in \R^{n}_x\times M^k_y$.
\\
\\
{\em First step: $\exists \tilde \lambda>0 \hbox{ s.t. } 
\forall v\in {\mathcal M}_{n,M^k, \alpha, \lambda},
v(x, y)>0 \hbox{ we have }
\nabla_y v=0, \hbox{ } \forall \lambda>\tilde \lambda$}
\\
\\
Assume that the conclusion is false then 
there exists
$\lambda_j\rightarrow \infty$ such that
$u_{\lambda_j}(x, y)\in  M_{n,M^k, \alpha, \lambda_j}, u_{\lambda_j}(x, y)>0
\hbox{ and } \nabla_y u_{\lambda_j}\neq 0$.
This is absurd due to Lemma \ref{unod}.
\\
\\
{\em Second step: conclusion}
\\
\\
We define 
$$\lambda^*=\inf_\lambda \{\lambda>0| \nabla_y v=0 
\hbox{ } \forall v\in {\mathcal M}_{n,M^k, \alpha, \lambda}\}$$
By the first step $\lambda^*<\infty$.
Moreover it is easy to deduce that if $\lambda>\lambda^*$ then the minimizers 
of the problem $J_{n,M^k, \alpha, \lambda}$
are precisely the same minimizers of the problem
$I_{n, \alpha}^{1/\sqrt {vol(M^k)}}$, which 
in turn are characterized in section \ref{cazlions} (hence we get \eqref{prop12}).
\\

Next we prove that $\lambda^*>0$. It is sufficient to show that 
\begin{equation}\label{lblambstar}
\lim_{\lambda\rightarrow 0} J_{n, M^k, \alpha, \lambda}< vol (M^k)I_{n,\alpha}^{1/\sqrt{ vol (M^k)}}
\end{equation}
(see \eqref{Id-1} and \eqref{JnMk} for a definition of the quantities involved in the inequality above).
Let us fix $\rho(y)\in C^\infty (M^k)$ such that $$\int_{M^k} |\rho|^2 dvol_{M^k_y}=1$$
and $\rho^2(y_0)\neq \frac 1{vol(M^k)}$ for some $y_0\in M^k$
(i.e. $\rho(y)$ is not identically constant).
Then we introduce the functions
$$\psi(x, y)=\rho(y)^{4/(4-\alpha n)}Q (\rho(y)^{(2 \alpha)/(4-\alpha n)} x)$$
where $Q(x)$ is the unique radially symmetric minimizer for
$I_{n,\alpha}^{1/\sqrt{vol(M^k)}}$.
Then we get
$$\|\psi(x, y)\|_{L^2_x}^2=(\rho(y))^2
\hbox{ and } {\mathcal E}_{n, \alpha} (\psi(x, y))
=I^{1}_{n,\alpha} (\rho(y))^{\frac{8+4\alpha-2\alpha n}{4-\alpha n}}
$$
and as a consequence we deduce
$$\int_{M^k_y} \int_{\R^n_x} \Big (\frac 12 |\nabla_x \psi(x, y)|^2 - 
\frac 1{2+\alpha}|\psi(x, y)|^{2+\alpha} \Big ) dx dvol_{M^k_y}$$
$$=I^{1}_{n,\alpha} \int_{M^k_y}(\rho(y))^{\frac{8+4\alpha-2\alpha n}{4-\alpha n}} d vol_{M^k_y}$$
$$<I^{1}_{n,\alpha} \Big (\int_{M^k} (\rho(y))^2 dvol_{M^k_y}\Big )^\frac{4-\alpha n +2\alpha}{4-\alpha n}
vol (M^k)^{-\frac{2\alpha}{4-\alpha n}}=I^{1}_{n,\alpha} vol (M^k)^{-\frac{2\alpha}{4-\alpha n}}$$
where at the last inequality we have used
the fact that $I^{1}_{n,\alpha}<0$ in conjunction with the H\"older inequality 
(moreover we get the inequality $<$ since by hypothesis $\rho(y)$ is not identically constant).
As a byproduct we get
$$\lim_{\lambda\rightarrow 0}{\mathcal E}_{n, M^k, \alpha, \lambda} (\psi(x, y))
<I^{1}_{n,\alpha} vol (M^k)^{-\frac{2\alpha}{4-\alpha n}}= 
vol (M^k)I_{n,\alpha}^{1/\sqrt{ vol (M^k)}}$$
(where we have used \eqref{potespl}) which in turn implies \eqref{lblambstar}.
\\

Let us finally prove \eqref{prop14}. It is sufficient to show that 
if $v\in {\mathcal M}_{n,M^k, \alpha, \lambda}$ for $\lambda<\lambda^*$ then $\nabla_y v\neq 0$.
Assume by the absurd that it is false, then we get
$\lambda_1<\lambda^*$ and 
$v_1\in {\mathcal M}_{n,M^k, \alpha, \lambda_1}$ such that $\nabla_y v_1=0$.
Arguing as above it implies that 
\begin{equation}\label{dabbic}J_{n,M^k, \alpha, \lambda_1}=vol(M^k)
I_{n, \alpha}^{1/\sqrt {vol(M^k)}}
\end{equation}
On the other hand by definition of $\lambda^*$ there exists 
$\lambda_2\in (\lambda_1, \lambda^*]$
and $v_2 \in {\mathcal M}_{n,M^k, \alpha, \lambda_2}$
such that $\nabla_y v_2\neq 0$.
As a consequence we deduce that
$$J_{n,M^k, \alpha, \lambda_1}< {\mathcal E}_{n, M^k, \alpha, \lambda_2} (v_2) 
=J_{n,M^k, \alpha, \lambda_2}\leq vol(M^k)I_{n, \alpha}^{1/\sqrt {vol(M^k)}}
$$
where at the last step we have used \eqref{minim}.
Hence we get a contradiction with \eqref{dabbic}.

\hfill$\Box$

\section{Proof of theorem \ref{variational}}\label{mainvar}

In the sequel the homogeneity
of the euclidean space $\R^n$ will play a key role. 
Due to this property we shall be able to reduce the proof of Theorem
\ref{variational} to the problem studied in the previous section.\\
In view of section \ref{cazlions}
it is sufficient to prove that there exists $\rho^*>0$ such that
\begin{equation}\label{implica}v\in {\mathcal M}_{n, M^k, \alpha}^\rho \hbox{ implies }
\nabla_y v=0 \hbox{ for } \rho<\rho^* 
\end{equation}
and
\begin{equation}\label{implica2}
v\in {\mathcal M}_{n, M^k, \alpha}^\rho \hbox{ implies }
\nabla_y v\neq 0 \hbox{ for } \rho>\rho^* 
\end{equation}
By an elementary computation we have that
the map
$$S_{1} \ni u\rightarrow \rho^{4/(4-\alpha n)} u (\rho^{2\alpha /(4-\alpha n)} x, y)\in 
S_\rho$$
where $$S_{\lambda}=\{v\in H^1(\R^n\times M^k)| \|v\|_{L_{x,y}^2}=\lambda\}$$ 
is a bijection.
Moreover we have
$${\mathcal E}_{n,M_k, \alpha}(\rho^{4/(4-\alpha n)} u (\rho^{2\alpha /(4-\alpha n)} x, y))=
\rho^{(8-2\alpha n)/(4-\alpha n)}\int_{M^k_y} \int_{\R^{n}_x} |\nabla_y u|^2 dxdvol_{M_y^k}$$$$+
\rho^{(8-2\alpha n+4\alpha)/(4-\alpha n)} \int_{M^k_y} \int_{\R^{n}_x} |\nabla_x u|^2 dxdvol_{M_y^k} 
$$$$- \rho^{(8-2\alpha n+4\alpha)/(4-\alpha n)}\frac{1}{2+\alpha} \int_{M_y^k} 
\int_{\R^n_x} |u|^{2+\alpha} dxdvol_{M_y^k}
$$
$$= \rho^{(8-2\alpha n+4\alpha)/(4-\alpha n)} \Big ( 
\frac 12 \rho^{-4\alpha/(4-\alpha n)} 
\int_{M^k_y} \int_{\R^{n}_x} |\nabla_y u|^2 dxdvol_{M_y^k}$$$$+
\frac 12 \int_{M^k_y} \int_{\R^{n}_x} |\nabla_x u|^2 
- \frac{1}{2+\alpha} |u|^{2+\frac 4d} dxdvol_{M_y^k}
\Big )$$
In particular \eqref{implica} 
and \eqref{implica2} are satisfied provided that there exists $\rho^*>0$ such that
\begin{equation}\label{implicap}v\in {\mathcal M}_{n, M^k, \alpha, \rho^{-4\alpha/(4-\alpha n)}} \hbox{ implies }
\nabla_y v=0 \hbox{ for } \rho<\rho^* 
\end{equation}
and
\begin{equation}\label{implica2p}
v\in {\mathcal M}_{n, M^k, \alpha, \rho^{-4\alpha/(4-\alpha n)}} \hbox{ implies }
\nabla_y v\neq 0 \hbox{ for } \rho>\rho^* 
\end{equation}
that in turn follow by Theorem \ref{rigidity}.

\section{Proof of Theorem \ref{cauchy}}

The main tool we use is the following Strichartz 
type estimates (whose proof follows by \cite{TV}).
\begin{prop}\label{Strich}
For every manifold $M^k_y$, $n\geq 1$ and $p, q\in [2, \infty]$
such that:
$$\frac 2p + \frac{n}q=\frac {n}2, \hbox{ } (p, n)\neq (2, 2)$$ 
there exists $C>0$ such that
\begin{equation}\label{firststri}\|e^{-{i}t\Delta_{x,y}} f
\|_{L^p_tL^q_xH^1_y} +
\Big \|\int_0^t e^{-{i}(t-s)\Delta_{x,y}} F(s) ds\Big \|_{L^p_tL^q_x H^1_y}
\end{equation}$$\leq C \Big (\|f \|_{L^2_{x}H^1_y}+
\|F\|_{L^{p'}_tL^{q'}_xH^1_y}\Big );$$
\begin{equation}\label{secondstri}\|\nabla_x 
e^{-{i}t\Delta_{x,y}} f\|_{L^p_tL^q_xL^2_y} +
\Big \| \nabla_x \int_0^t e^{-{i}(t-s)\Delta_{x,y}} F(s) ds
\Big \|_{L^p_tL^q_x L^2_y}
\end{equation}
$$\leq C \Big (\|\nabla_x f \|_{L^2_{x}L^2_y}+
\|\nabla_x F\|_{L^{p'}_tL^{q'}_xL^2_y}\Big )$$
and
\begin{equation}\label{secthndstri}
\Big \|\int_0^t e^{-{i}(t-s)\Delta_{x,y}} F(s) ds\Big \|_{L^p_tL^q_x L^2_y}
\leq C \|F\|_{L^{p'}_tL^{q'}_xL^2_y}\end{equation}
Moreover 
\begin{equation}\label{stribis}\|e^{-{i}t\Delta_{x,y}} 
f\|_{L^\infty_tL^2_xH^1_y} +
\Big \|\int_0^t e^{-{i}(t-s)\Delta_{x,y}} F(s) ds
\Big \|_{L^\infty_tL^2_x H^1_y}
\end{equation}$$\leq C \Big (\|f \|_{L^2_{x}H^1_y}
+\|F\|_{L^{p'}_tL^{q'}_xH^1_y}\Big )$$
and
\begin{equation}\label{stribis2}
\|\nabla_x e^{-{i}t\Delta_{x,y}} f\|_{L^\infty_tL^2_xL^2_y} +
\Big \| \nabla_x \int_0^t e^{-{i}(t-s)\Delta_{x,y}} F(s) ds
\Big \|_{L^\infty_tL^2_x L^2_y}
\end{equation}
$$\leq C \Big (\|\nabla_x f \|_{L^2_{x}L^2_y}+
\|\nabla_x F\|_{L^{p'}_tL^{q'}_xL^2_y}\Big )$$
\end{prop}
Next we shall use the norms  
$\|.\|_{X_T}$ and $\|.\|_{Y_T}$ introduced in \eqref{XT} and \eqref{YT}
for time dependent functions.
We also introduce the space $Z_T$ whose norm is defined by
$$\| v\|_{Z_T}\equiv\|v\|_{X_T} + \|v\|_{Y_T}$$
and the nonlinear operator
associated to the Cauchy problem
\eqref{NLS}:
$${\mathcal T}_{\varphi}(u)\equiv
e^{-{i}t\Delta_{x,y}} \varphi + \int_0^t e^{-{i}(t-s)\Delta_{x,y}} u(s)|u(s)|^{\alpha} ds$$
We split the proof of Theorem \ref{cauchy}
in several steps.

\subsection{Local Well Posedness}
This subsection is devoted to the proof of the following fact:
$$\forall \varphi\in H^{1}(\R^n \times M^1) \hbox{  } \exists T=T(\|\varphi\|_{H^1(\R^n \times M^1)})>0
$$$$\hbox{ and } \exists! v(t,x)\in Z_{T}\cap {\mathcal C}((-T, T);H^1(\R^n \times M^1))
$$$$\hbox{ s.t. }
{\mathcal T}_{\varphi} v(t)=v(t)\ \hbox{ } \forall t\in (-T, T)$$
\\
{\em First step: 
$$\forall \varphi\in H^{1}(\R^n \times M^1) \hbox{  } \exists T=T(\|\varphi\|_{H^1(\R^n \times M^1)})>0, 
R=R(\|\varphi\|_{H^1(\R^n \times M^1)})>0
\hbox{ s.t. }$$$$ {\mathcal T}_\varphi(B_{Z_{\tilde T}}(0, R))
\subset B_{Z_{\tilde T}}(0, R)
\hbox{ } \forall \tilde T< T$$
}
\\
\\
First we estimate the nonlinear term:
$$\|u|u|^{\alpha}\|_{{L^{p'}_t L^{q'}_x H^1_y}}
\leq \left \| \| u^{\alpha}(t,x,.)\|_{L^\infty_y} 
\|u(t,x,.) \|_{H^1_y} \right \|_{L^{p'}_t L^{q'}_x }
$$
(where $(p,q)$ is the couple in \eqref{XT}
and \eqref{YT}) and after application of the H\"older inequality in $(t, x)$ we get
$$...\leq  \|u\|_{L^{p}_tL^{q}_xH^1_y}
\|u\|_{L^{\alpha \tilde p}_t L^{\alpha \tilde q}_xL^\infty_y}^{\alpha}
$$$$\leq C \|u\|_{L^{p}_tL^{q}_xH^1_y}
\|u\|_{L^{\alpha \tilde p}_t L^{\alpha \tilde q}_x H^1_y}^{\alpha}
$$
where we have used the embedding
$H^1_y\subset L^\infty_y$ and we have chosen
$$\frac 1{\tilde p}+ \frac 1{p}= 1-\frac{1}{p}$$
$$\frac 1{\tilde q}+ \frac 1{q}= 1-\frac{1}{q}$$
By direct computation we have:
\begin{equation}\label{monoton}\alpha \tilde q= q \hbox{ and }
\alpha  \tilde p<p\end{equation} 
By combining the nonlinear estimate above with \eqref{firststri}, \eqref{monoton} and
the H\"older inequality 
(in the time variable) we get: 
\begin{equation}\label{icsd}
\|{\mathcal T}_\varphi u\|_{X_{T}}
\leq  C (\|\varphi\|_{L^2_xH^1_y} + 
T^{a(d)} \|u\|_{X_{T}}^{1+\alpha})\end{equation}
with $a(d)>0$.\\
Arguing as above get
$$\|\nabla_x (u|u|^{\alpha}) \|_{L^{p'}_t L^{q'}_x L^2_y}
\leq C \|\nabla_x u\|_{L^{p}_tL^{q}_x L^2_y} 
\|u^{\alpha}\|_{L^{\tilde p}_t L^{\tilde q}_xL^\infty_y}
$$
$$\leq C \|u\|_{Y_{T}} \|u\|_{L^{\alpha \tilde p}_t 
L^{\alpha \tilde q}_x H^1_y }^{\alpha}$$
where $\tilde p$ and $\tilde q$ are as above and we have used the embedding
$H^1_y \subset L^{\infty}_y$.
As a consequence of this estimate and \eqref{secondstri}
we get:
\begin{equation}\label{icsd2}
\|{\mathcal T}_\varphi u\|_{Y_{T}}
\leq  C (\|\nabla_x \varphi\|_{L^2_{x,y}} + 
T^{a(d)} \| u\|_{Y_{T}}\|u\|_{X_{T}}^{\alpha})\end{equation}
with $a(d)>0$.\\
By combining \eqref{icsd} with \eqref{icsd2}
we get
\begin{equation*}
\|{\mathcal T}_\varphi u\|_{Z_{T}}
\leq  C (\|\varphi\|_{H^1(\R^n\times M^1)} + 
T^{a(d)} \| u\|_{Z_{T}}\|u\|_{Z_{T}}^{\alpha})\end{equation*}
The proof follows by a standard continuity argument.
\\
\\
Next we introduce the norm
$$\|w(t,x,y) \|_{\tilde Z_T}\equiv \|w(t,x,y)\|_{L^{p}((-T,T); L^{q}_x L^2_y)}$$
and we shall prove the following.
\\
\\
{\em Second step: let $T, R>0$ as in the previous step
then $$\exists T'=T'(\|\varphi\|_{H^1(\R^n \times M^1)})<T
\hbox{ s.t. }  {\mathcal T}_\varphi$$
$$\hbox{ is a contraction on } B_{Z_{T'}}(0, R) \hbox{ endowed with the norm } 
\|.\|_{\tilde Z_{T'}}$$}
\\
\\
It is sufficient to prove:
\begin{equation}\label{icsdy}
\|{\mathcal T}_\varphi v_1 - {\mathcal T}_\varphi v_2\|_{\tilde Z_{T}}
\leq  C 
T^{a(d)} \|v_1-v_2\|_{\tilde Z_{T}} \sup_{i=1,2} \{ \|v_i\|_{Z_{T}}\}^{\alpha}
\end{equation}
with $a(d)>0$. Notice that we have
$$\|v_1|v_1|^\alpha-v_2|v_2|^\alpha\|_{L^{p'}((-T, T); L^{q'}_x L^2_y)}
$$$$\leq C \Big \|\|v_1-v_2\|_{L^2_y} (\|v_1\|_{L^\infty_y} 
+ \|v_2\|_{L^\infty_y})^\alpha\Big\|_{L^{p'}((-T,T); L^{q'}_x)}
$$$$\leq C 
T^{a(d)} \|v_1-v_2\|_{\tilde Z_{T}} \sup_{i=1,2} \{ \|v_i\|_{Z_{T}}\}^{\alpha} 
$$
where we have used the Sobolev embedding $H^1_y\subset L^\infty_y$
and the H\"older inequality in the same spirit as in
the proof of \eqref{icsd} and \eqref{icsd2}.
We conclude by combining the estimate above with the Strichartz estimate
\eqref{secthndstri}.
\\
\\
{\em Third step: existence and uniqueness of solution in $Z_{T'}$ where $T'$ 
is as in the previous step}
\\
\\
We apply the contraction principle to the map
${\mathcal T}_\varphi$ defined on the complete space
$B_{Z_{T'}}(0, R)$ endowed with the topology induced by $\|.\|_{\tilde Z_{T'}}$.
It is well-known that this space is complete.
\\
\\
{\em Fourth step: regularity of the solution}
\\
\\
By combining the previous steps with the fixed point argument we get
the existence of a solution $v\in Z_{T'}$. In order to get the regularity
$v\in {\mathcal C}((-T', T'); H^1(\R^n \times M^1))$
it is sufficient to argue as in the first step
(to estimate the nonlinearity) in conjugation with the  
Strichartz estimates \eqref{stribis} and \eqref{stribis2}.

\subsection{Global Well Posedness}

Next we prove that the local solution (whose existence has been proved above)
cannot blow-up in finite time.
The argument is standard and follows from the 
conservation laws:
\begin{equation}\label{coSE}
\|u(t)\|_{L^2_{x,y}}\equiv \|\varphi\|_{L^2_{x,y}}
\end{equation}
\begin{equation}\label{energyconservati}{\mathcal E}_{n, M^1, \alpha} (u(t)) 
+ \frac 12 \|u(t)\|_{L^2_{x,y}}^2\equiv 
{\mathcal E}_{n,M^1, \alpha}(\varphi) + \frac 12 \|\varphi\|_{L^2_{x,y}}^2
\end{equation}
where 
${\mathcal E}_{n,M^1, \alpha}$ is defined in \eqref{L2peri}.
By the Gagliardo Nirenberg inequality we deduce
$${\mathcal E}_{n, M^1, \alpha} (u(t)) 
+ \frac 12 \|u(t)\|_{L^2_{x,y}}^2\geq \frac 12 \|u(t)\|_{H^1(\R^n \times M^1)}^2
- C \|u(t)\|_{L^2_{x,y}}^{2+\alpha - \mu} \|u(t)\|_{H^1(\R^n \times M^1)}^\mu
$$
for a suitable $\mu\in (0, 2)$. By combining the estimate above with \eqref{coSE}
and \eqref{energyconservati} we get
$$\frac 12 \|u(t)\|_{H^1(\R^n \times M^1)}^2
- C \|\varphi\|_{L^2_{x,y}}^{2+\alpha - \mu} \|u(t)\|_{H^1(\R^n \times M^1)}^\mu
\leq 
{\mathcal E}_{n,M^1, \alpha}(\varphi) + \frac 12 \|\varphi\|_{L^2_{x,y}}^2
$$
Since $\mu\in (0, 2)$ it implies that $\|u(t)\|_{H^1(\R^n\times M^1)}$
cannot blow-up in finite time.

\section{Appendix}
For the sake of completeness we prove in this appendix Theorems \ref{variationalappendix}
and \ref{stabilitysubordinate}. Our argument is heavily inspired by the work \cite{CL}
even if, in our opinion, the following presentation of Theorem \ref{variationalappendix}
is simpler compared with the original one.\\
\\
{\bf Proof of Theorem \ref{variationalappendix}}
For any given $\rho>0$ we shall denote by 
$u_{j,\rho}\in H^1(\R^n \times M^k)$
any constrained minimizing sequence, i.e.:
\begin{equation}\|u_{j, \rho}\|_{L^2_{x,y}}=\rho
\hbox{ and } 
\lim_{j\rightarrow \infty }{\mathcal E}_{n, M^k, \alpha}(u_{j, \rho})=
K_{n, M^k, \alpha}^\rho
\end{equation}
Next we split the proof in many steps.
\\
\\
{\em First step:  $K^{\rho}_{n,M^k, \alpha}>-\infty$ and
$\sup_{j}\|u_{j, \rho}\|_{H^1_{x,y}}<\infty$, $\forall \rho>0$}\\
\\
\\
By the classical Gagliardo Nirenberg inequality (see \eqref{GagN}) we get
the existence of $\mu\in (0, 2)$ such that 
$$\mathcal E_{n, M^k, \alpha}(u_{j, \rho})+\frac 12 \rho^2$$$$
\geq \frac 12 \int_{M^k_y} \int_{\R^{n}}
(|\nabla_{x,y} u_{j, \rho}|^2 
+ |u_{j, \rho}|^2) dxdvol{M_y^k}
- C(\rho) \|u_{j, \rho}\|_{H^1(\R^m \times M^k)}^\mu
$$$$\geq \inf_{t>0} (1/2t^2 - C(\rho) t^\mu)>-\infty
$$
The conclusion follows by a standard argument.
\\
\\
{\em Second step: the map
$(0, \infty)\ni \rho\rightarrow K^\rho_{n, M^k, \alpha}$ is continuous}
\\
\\
Fix $\rho\in (0, \infty)$ and let $\rho_j \rightarrow \rho$.
Then we have
$$K^{\rho_j}_{n, M^k, \rho}\leq {\mathcal E}_{n, M^k, \alpha}\Big (
\frac{\rho_j}{\rho}u_{j, \rho}\Big )
=
$$
$$
\Big (\frac{\rho_j}{\rho}\Big )^2 \Big(
\frac 12 \|\nabla_{x,y} u_{j, \rho}\|_{L^2_{x,y}}^2 
- \frac {1}{2+\alpha} \Big ( \frac{\rho_j}{\rho} \Big)^{\alpha}
\|u_{j, \rho}\|_{L^{2+\alpha}_{x,y}
}^{2+\alpha} \Big)$$
$$=
\Big (\frac{\rho_j}{\rho}\Big )^2 \Big (\frac 12 
\|\nabla_{x,y} u_{j, \rho}\|_{L^2_{x,y}}^2 
- \frac 1{2+\alpha} \|u_{j, \rho}\|_{L^{2+\alpha}_{x,y}}^{2+\alpha} \Big)$$$$
+ \frac 1{2+\alpha} \Big(\frac{\rho_j}{\rho}\Big)^2
\Big (1- \left ( \frac{\rho_j}{\rho}\right)^\alpha\Big )
\|u_{j, \rho}\|_{L^{2+\alpha}_{x,y}}^{2+\alpha}
$$
$$= \Big (\frac 12 \|\nabla_{x,y} u_{j, \rho}\|_{L^2_{x,y}}^2 
- \frac 1{2+\alpha} \|u_{j, \rho}
\|_{L^{2+\alpha}_{x,y}}^{2+\alpha}\Big )$$$$+
\Big (\Big (\frac{\rho_j}{\rho}\Big )^2 -1\Big )
\Big (\frac 12 \|\nabla_{x,y} u_{j, \rho}\|_{L^2_{x,y}}^2
- \frac 1{2+\alpha} \|u_{j, \rho}\|_{L^{2+\alpha}_{x,y}}^{2+\alpha}\Big )$$$$
+ \frac 1{2+\alpha} \Big(\frac{\rho_j}{\rho}\Big)^2
\Big (1- \Big ( \frac{\rho_j}{\rho}\Big )^\alpha \Big )
\|u_{j, \rho}\|_{L^{2+\alpha}_{x,y}}^{2+\alpha}
$$
Since we are assuming that $\rho_j \rightarrow \rho$
and $\sup_n \|u_{j, \rho}\|_{H^1(\R^n \times M^k)}<\infty$ 
(see the first step)
we get
$$\limsup_{j\rightarrow \infty} K_{n,M^k, \alpha}^{\rho_j}
\leq K_{n,M^k, \alpha}^{\rho}$$
To prove the opposite inequality let us fix
$u_j\in H^1(\R^n \times M^k)$ such that\begin{equation}\label{choice}\|u_j\|_{L^2_{x,y}}
=\rho_j
\hbox{ and } {\mathcal E}_{n,M^k, \alpha}(u_j)
<K^{\rho_j}_{n,M^k, \alpha}+\frac 1j\end{equation}
By looking at the proof of the first step we also deduce that
$u_j$ can be chosen in such a way that 
\begin{equation}\label{boundunif}
\sup_{j} \|u_j\|_{H^1(\R^n\times M^k)}<\infty
\end{equation}
Then we can argue as above an we get
$$K^{\rho}_{n,M^k, \alpha}\leq {\mathcal E}_{n,M^k, \alpha}
\left (\frac{\rho}{\rho_j}u_{j}\right )
$$
$$= \Big (\frac 12 \|\nabla_{x,y} u_{j}\|_{L^2_{x,y}}^2 
- \frac 1{2+\alpha} \|u_{j}\|_{L^{2+\alpha}_{x,y}}^{2+\alpha}\Big )$$$$+
\Big (\Big (\frac{\rho}{\rho_j}\Big )^2 -1\Big )
\Big (\frac 12 \|\nabla_{x,y} u_{j}\|_{L^2_{x,y}}^2
- \frac 1{2+\alpha} \|u_{j}\|_{L^{2+\alpha}_{x,y}}^{2+\alpha}\Big )$$$$
+ \frac 1{2+\alpha} \Big(\frac{\rho}{\rho_j}\Big)^2
\Big (1- \left ( \frac{\rho}{\rho_j}\right)^{\alpha}\Big )
\|u_{j}\|_{L^{2+\alpha}_{x,y}}^{2+\alpha}
$$
By using \eqref{choice}, \eqref{boundunif} 
and the assumption $\rho_j\rightarrow \rho$ we get
$$K_{n,M^k, \alpha}^{\rho}\leq \liminf_{j\rightarrow \infty} 
K_{n,M^k, \alpha}^{\rho_j}$$
\\
\\
{\em Third step: for every $\rho>0$ we have (up to subsequence)
$\inf_j \|u_{j, \rho}\|_{L^{2+\alpha}_{x,y}}>0$}
\\
\\
It is sufficient to prove that
$K_{n,M^k, \alpha}^{\rho}<0.$ In fact we have
\begin{equation}\label{chain}
K_{n,M^k, \alpha}^{\rho}\leq vol(M^k)
{\mathcal E}_{n,\alpha}(u_{n, \omega, \alpha})=vol(M^k) 
I_{n,\alpha}^{\rho/\sqrt {vol(M^k)}}<0
\end{equation}
where ${\mathcal E}_{n,\alpha}$ is the energy defined in \eqref{energyd-1}
and $\omega$ is chosen in such a way that $\|u_{n,\omega, \alpha}\|_{L^2_x}=
\frac{\rho}{\sqrt {vol(M^k)}}$.
Notice that in \eqref{chain}
we have used \eqref{sign} and \eqref{idfond}.
\\
\\
{\em Fourth step: for any minimizing sequence $u_{j, \rho}$ there exists
$\tau_j\in \R^{n}$ s.t. (up to subsequence)
$u_{j, \rho}(x+\tau_j,y)$ has a weak limit $\bar u\neq 0$}
\\
\\
We have the following
localized Gagliardo Nirenberg inequality:
\begin{equation}\label{GNprecised}
\|v\|_{L^{2+4/(n+k)}_{x,y}}\leq C \sup_{x\in \R^{n}} 
\Big (\|v\|_{L^2_{Q^{n}_x\times M^k}}\Big )^{2/(n+k+2)} 
\|v\|_{H^1(\R^n\times M^k)}^{(n+k)/(n+k+2)}
\end{equation}
where $$Q^{n}_x=x+ [0,1]^{n} \hbox{ } \forall x\in \R^{n}$$
The estimate above can be proved as follows (see \cite{L}
for a similar argument 
on the flat space $\R^{d+k}$).
We fix $x_h\in \R^{n}$ in such a way that
$\bigcup_h Q^{n}_{x_h}=\R^{n}$ and
$meas_{n} (Q^{n}_{x_i}\cap Q^{n}_{x_j})=0$ for $i\neq j$
where $meas_{n}$ denotes the Lebesgue measure in $\R^{n}$. 
By the classical Gagliardo Nirenberg inequality we get:
$$
\|v\|_{L^{2+4/(n+k)}_{Q^{n}_{x_h}\times M^k}}^{2+4/(n+k)}
\leq C \|v\|_{L^{2}_{Q^{n}_{x_h}\times M^k}}^{4/(n+k)}
\|v\|_{H^1(Q^{n}_{x_h}\times M^k)}^2$$
The proof of \eqref{GNprecised}
follows by taking the sum of the previous estimates 
on $h\in \N$.\\
Due to the boundedness of $u_{j,\rho}$ in $H^1(\R^m \times M^k)$
(see the first step) we deduce by \eqref{GNprecised}
that
\begin{equation}\label{inf2d}
0<\epsilon_0=\inf_j \|u_{j, \rho}\|_{L^{2+4/(n+k)}_{x,y}}\leq C \sup_{x\in \R^{n}} 
\|u_{j, \rho}\|_{L^2_{Q^{n}_x
\times M^k}}^{2/(n+k+2)}
\end{equation}
(the l.h.s. above follows by combining the H\"older inequality with the third step).
The proof can be concluded by the Rellich compactness theorem
once we choose a
sequence $\tau_j\in \R^{n}_x$ in such a way that
$$\inf_j \|u_{j, \rho}\|_{L^2_{Q^{n}_{\tau_j}\times M^k}}> 0$$
(the existence of such a sequence $\tau_j$ follows by 
\eqref{inf2d}).
\\
\\
{\em Fifth step: the map
$(0, \bar \rho) \ni \rho \rightarrow \rho^{-2} K_{n,M^k, \alpha}^{\rho}$ 
is strictly decreasing}
\\
\\
Let us fix $\rho_1<\rho_2$ 
and $u_{j, \rho_1}$ a minimizing sequence for
$K_{n,M^k, \alpha}^{\rho_1}$. Then we have
$$K_{n,M^k, \alpha}^{\rho_2}\leq 
{\mathcal E}_{n,M^k, \alpha}\Big (\frac{\rho_2}{\rho_1} u_{j, \rho_1}\Big )
$$$$= \Big (\frac{\rho_2}{\rho_1}\Big )^2 \Big (\frac 12 
\|\nabla_{x,y} u_{j, \rho_1}\|_{L^2_{x,y}}^2 
- \frac {1}{2+\alpha} \Big ( \frac{\rho_2}{\rho_1}\Big)^{\alpha}
\|u_{j, \rho_1}\|_{L_{x,y}^{2+\alpha}
}^{2+\alpha} \Big)$$
$$= \Big (\frac{\rho_2}{\rho_1}\Big )^2 \Big 
(\frac 12 \|\nabla_{x,y} u_{j, \rho_1}\|_{L^2_{x,y}}^2 
- \frac 1{2+\alpha} \|u_{j, \rho_1}\|_{L^{2+\alpha}_{x,y}}^{2+\alpha}  \Big)$$$$
+ \frac 1{2+\alpha} \Big(\frac{\rho_2}{\rho_1}\Big)^2
\Big (1- \Big ( \frac{\rho_2}{\rho_1}\Big)^{\alpha}\Big ) \|u_{j, \rho_1}
\|_{L^{2+\alpha}_{x,y}}^{2+\alpha}$$$$ \leq 
\Big (\frac{\rho_2}{\rho_1}\Big )^2 \Big (\frac 12 
\|\nabla_{x,y} u_{j, \rho_1}\|_{L^2_{x,y}}^2 
- \frac 1{2+\alpha} \|u_{j, \rho_1}\|_{L^{2+\alpha}_{x,y}}^{2+\alpha} \Big)$$$$
+ \frac 1{2+\alpha} \Big(\frac{\rho_2}{\rho_1}\Big)^2
\Big (1- \Big ( \frac{\rho_2}{\rho_1}\Big)^{\alpha}\Big )
\inf_j \|u_{j, \rho_1}\|_{L^{2+\alpha}_{x,y}}^{2+\alpha}
$$
By recalling (see the third step) that
$\inf_j \|u_{j, \rho_1}\|_{L^{2+\alpha}_{x,y}}^{2+\alpha}>0$ we get
$$K_{n,M^k, \alpha}^{\rho_2}<\Big 
(\frac{\rho_2}{\rho_1}\Big )^2 K^{\rho_1}_{n,M^k, \alpha}$$
\\
\\
{\em Sixth step: let 
$\bar u$ be as in the fourth step, then $\|\bar u\|_{L^2_{x,y}}=\rho$}
\\
\\
Up to subsequence we get:
$$u_{j,\rho}(x+\tau_{j},y)\rightarrow \bar u(x,y)\neq 0 
\hbox{ a.e. } (x,y)\in \R_x^{n}\times M^k_y$$
and hence by the Brezis-Lieb lemma (see \cite{BL}) we get
\begin{equation}\label{BL}
\|u_{j, \rho}(x+\tau_{j},y)-\bar u(x,y)\|_{L^{2+\alpha}_{x,y}}^{2+\alpha}
\end{equation}$$=\|u_{j, \rho}(x+\tau_{j},y)\|_{L^{2+\alpha}_{x,y}}^{2+\alpha} - 
\|\bar u(x,y)\|^{2+\alpha}_{L^{2+\alpha}_{x,y}}+ o(1)$$
Assume that $\|\bar u\|_{L^2_{x,y}}=\theta$,
our aim is to prove $\theta=\rho$. Since $\bar u\neq 0$ necessarily $\theta>0$.
Notice that since $L^2_{x,y}$ is an Hilbert space we have
\begin{equation}\label{pger}\rho^2=\|u_{j, \rho}(x+\tau_{j},y)
\|_{L^2_{x,y}}^2\end{equation}$$=\|u_{j, \rho}(x+\tau_{j},y)-\bar u(x,y)\|_{L^2_{x,y}}^2 
+\|\bar u(x,y)\|_{L^2_{x,y}}^2 + o(1)$$
and hence
\begin{equation}\label{orthgo}
\|u_{j, \rho}(x+\tau_{j},y)-\bar u(x,y)\|_{L^2_{x,y}}^2 =\rho^2 - \theta^2+ o(1)
\end{equation}
By a similar argument
\begin{equation}\label{orthgo2}
\int_{M^k_y} \int_{\R_x^{n}}|\nabla_{x} (u_{j, \rho}(x+\tau_{j},y)) 
- \nabla_{x} \bar u(x,y)|^2 dxdy \end{equation}
$$
+ \int_{M^k_y} \int_{\R_x^{n}} |\nabla_{y} (u_{j, \rho}(x+\tau_{j},y)) 
- \nabla_{y} \bar u(x,y)|^2 dxdvol_{M^k_y}
$$$$
+ \int_{M_y^k} \int_{\R_x^{n}} 
(|\nabla_{x} \bar u(x,y)|^2 + |\nabla_y \bar u(x,y)|^2) dxdvol_{M^k_y}
$$
$$=
\int_{M^k_y} \int_{\R_x^{n}}  (|\nabla_{x} (u_{j, \rho}(x+\tau_{j},y)|^2 +
|\nabla_y u_{j, \rho}(x+\tau_j,y)|^2)dxdvol_{M^k_y}+ o(1)
$$
By combining \eqref{orthgo2} with \eqref{BL} we get:
\begin{equation}\label{benic}
K_{n,M^k, \alpha}^{\rho}=\lim_{j\rightarrow \infty} {\mathcal E}_{n, M^k, \alpha}(u_{j, \rho}(x+\tau_{j},y))=
\end{equation}
$$\lim_{j\rightarrow \infty} {\mathcal E}_{n,M^k, \alpha}(u_{j, \rho}(x+\tau_{j},y) - \bar u(x,y)) 
+ {\mathcal E}_{n,M^k, \alpha}(\bar u)
$$
and we can continue the estimate as follows
$$...\geq K_{n,M^k, \alpha}^{\sqrt{\rho^2- \theta^2} + o(1)}+K_{n,M^k, \alpha}^{\theta}
$$
where we have used \eqref{orthgo}. Hence by using the second step 
we get $$K_{n,M^k, \alpha}^\rho\geq  
K_{n,M^k, \alpha}^{\sqrt{\rho^2- \theta^2}}+K_{n, M^k, \alpha}^{\theta}
$$
Assume that $\theta<\rho$, then by using the monotonicity proved in fifth step we get
$$
K_{n,M^k, \alpha}^\rho> \frac{\rho^2-\theta^2}{\rho^2} K_{n,M^k, \alpha}^\rho + \frac{\theta^2}{\rho^2} 
K_{n, M^k, \alpha}^\rho=K_{n,M^k, \alpha}^\rho$$
and we have an absurd.

\hfill$\Box$
\\
\\
{\bf Proof of Theorem \ref{stabilitysubordinate}}
Assume by the absurd that the conclusion is false, then there exists $\rho$
and two sequences $\varphi_j\in H^1(\R^n\times M^k)$
and $t_j\in \R$ such that
\begin{equation}\label{appeh}\lim_{j\rightarrow \infty} dist_{H^1(\R^n\times M^k)}( \varphi_j, \mathcal M_{n, M^k, \alpha}^\rho)=0
\end{equation}
and
\begin{equation}\label{appef}\liminf_{j\rightarrow \infty} 
dist_{H^1(\R^n\times M^k)}( u_{\varphi_j}(t_j), \mathcal M_{n, M^k, \alpha}^\rho)>0
\end{equation}
where $u_{\varphi_j}$ is the solution to \eqref{NLS} with Cauchy data $\varphi_j$.
By \eqref{appeh} we deduce the following informations:
$$\lim_{j\rightarrow \infty} \|\varphi_j\|_{L^2_{x,y}}=\rho
\hbox{ and } \lim_{j\rightarrow \infty}
{\mathcal E}_{n, M^k, \alpha}(\varphi_j)=K_{n,M^k, \alpha}^\rho$$
and hence due to the conservation laws satisfied by solutions to \eqref{NLS}
we get
$$\lim_{j\rightarrow \infty} \|u_{\varphi_j}(t_j)\|_{L^2_{x,y}}=\rho
\hbox{ and } \lim_{j\rightarrow \infty}
{\mathcal E}_{n, M^k, \alpha}(u_{\varphi_j}(t_j))=K_{n,M^k, \alpha}^\rho$$
In turn by an elementary computation we get:
$$\|\tilde u_j\|_{L^2_{x,y}}=\rho \hbox{ and } 
\lim_{j\rightarrow \infty}
{\mathcal E}_{n, M^k, \alpha}(\tilde u_j)=K_{n,M^k, \alpha}^\rho$$
(more precisely $\tilde u_j$ is constrained minimizing sequence for $K^\rho_{n,M^k, \alpha}$)
where
$$\tilde u_{j}=\rho \frac{u_{\varphi_j}(t_j)}{\|u_{\varphi_j}(t_j)\|_{L^2_{x,y}}}$$
Moreover by \eqref{appef} it is easy to deduce
$$\liminf_{j\rightarrow \infty} dist_{H^1(\R^n\times M^k)}( \tilde u_{j}, \mathcal M_{n, M^k, \alpha}^\rho)>0
$$
and it is in contradiction 
with the compactness of minimizing sequences for
$K^\rho_{n,M^k, \alpha}$ stated in Theorem \ref{variationalappendix}.

\hfill$\Box$

\bibliographystyle{plain}

\begin{thebibliography}{1}

\bibitem{BL} H. Brezis, E,  Lieb, {\it  A relation between pointwise convergence of func-tions and convergence of functionals}, 
Proc. AMS,  88, (1983) 486-490.
%
\bibitem{BGT1} N. Burq, P. Gerard, N. Tzvetkov,
{\it Strichartz inequalities and the nonlinear Schrodinger equation on compact manifolds}, Amer. J. Math. 126 (2004), 569--605. 
%
\bibitem{BGT2} N. Burq, P. Gerard, N. Tzvetkov,
{\it The Cauchy problem for the non linear Schrodinger equation on a compact manifold}, J. Nonlinear Math. Phys. 10 (2003), 12--27.
%
\bibitem{C} T. Cazenave, {\it Semilinear Schr\"odinger equations}, Courant Lecture Notes 10 (2003).
%
\bibitem{CL} T. Cazenave, P.L. Lions, {\it Orbital stability of standing waves for some nonlinear Schr\"odinger equations}, Comm. Math. Phys. 85 (1982) 549-561.
%
\bibitem{HTT}
S. Herr, D. Tataru, N. Tzvetkov,
{\it Strichartz estimates for partially periodic solutions to Schršdinger equations in 4d and applications}, preprint, http://arxiv.org/abs/1011.0591 
%
\bibitem{IP} A. Ionescu, B. Pausader,
{\it Global well-posedness of the energy-critical defocusing NLS on $\mathbb{R}\times\mathbb{T}^3$}, to appear in Comm. Math. Phys.
%
\bibitem{L} P.L. Lions,  {\it The concentration compactness principal in the calculus of variations. I and II.}, Ann. IHP, Analyse non lin\'eaire 1 (1984) 109-145, 223-283.
%
\bibitem{RT}  F. Rousset, N. Tzvetkov,
{\it Stability and instability of the KdV solitary wave under the KP-I flow}, to appear in Comm. Math. Phys.
%
\bibitem{Sogge} C. Sogge,
{\it Fourier integrals in classical analysis}, Cambridge University Press, 1993.
%
\bibitem{Ts} Y. Tsutsumi, {\it $L^2$ solutions for nonlinear Schr\"odinger equations and nonlinear groups}, Funkcial. Ekvac. 30 (1987) 115-125.
%
\bibitem{TV} N. Tzvetkov, N. Visciglia, {\it Small data scattering for the nonlinear Schr\"odinger equation on product spaces}, Comm. PDE 37 (2012) 125-135.
%
\bibitem{W} M. Weinstein, {\it Lyapunov stability of ground states of nonlinear dispersive evolution equations}, Comm. Pure. Appl. Maths. 29 (1986) 51-68. 

\end{thebibliography}

\end{document}